\documentclass[10pt,a4paper]{amsart}
\usepackage{amsmath}
\usepackage{amssymb}
\usepackage{epsf}
\usepackage{amscd}
\usepackage[all]{xy}
\usepackage{enumerate}
\usepackage{mathpazo,color}
\usepackage[normalem,normalbf]{ulem}
\usepackage[colorlinks=true,linkcolor=blue]{hyperref}
\usepackage{pdflscape,afterpage,capt-of,tabularx}

\allowdisplaybreaks[2]

\usepackage[hmargin=2cm,vmargin={1.5cm,2cm}]{geometry}

\numberwithin{equation}{section}

\newtheorem{theo}{Theorem}[section]

\newtheorem{prop}[theo]{Proposition}

\newtheorem{cor}[theo]{Corollary}

\newtheorem{lemma}[theo]{Lemma}

\newtheorem{dfn}[theo]{Definition}

\newtheorem{rquee}[theo]{{\itshape Remark}}
\newenvironment{rmk}{\begin{rquee} \normalfont}{\end{rquee}}

\newtheorem{pty}[theo]{Property}

\newtheorem{notationse}[theo]{{\itshape Notation}}
\newenvironment{nota}{\begin{notationse} \normalfont}{\end{notationse}}

\newcommand{\bd}{\begin{dfn}}
\newcommand{\ed}{\end{dfn}}
\newcommand{\bp}{\begin{prop}}
\newcommand{\ep}{\end{prop}}
\newcommand{\bt}{\begin{theo}}
\newcommand{\et}{\end{theo}}
\newcommand{\bc}{\begin{cor}}
\newcommand{\ec}{\end{cor}}
\newcommand{\bl}{\begin{lemma}}
\newcommand{\el}{\end{lemma}}
\newcommand{\br}{\begin{rmk}}
\newcommand{\er}{\end{rmk}}
\newcommand{\bpf}{\begin{proof}}
\newcommand{\epf}{\end{proof}}
\newcommand{\bn}{\begin{nota}}
\newcommand{\en}{\end{nota}}

\newcommand{\ppq}{\leqslant}
\newcommand{\pgq}{\geqslant}
\newcommand{\zz}{\mathbb{Z}}

\newcommand{\Bb}{\mathbb{B}}
\newcommand{\rad}{\operatorname{rad}\nolimits}
\newcommand{\soc}{\operatorname{soc}\nolimits}
\newcommand{\Top}{\operatorname{top}\nolimits}
\newcommand{\id}{\operatorname{id}\nolimits}
\newcommand{\Hom}{\operatorname{Hom}\nolimits}
\newcommand{\sHom}{\operatorname{\underline{Hom}}\nolimits}
\newcommand{\End}{\operatorname{End}\nolimits}

\newcommand{\smod}{\ensuremath{\operatorname{\underline{mod}}\nolimits}}
\newcommand{\Ext}{\operatorname{Ext}\nolimits}
\newcommand{\A}{\mathcal{A}}
\newcommand{\I}{\mathfrak{I}}
\newcommand{\lnd}{\Lambda_{n,d}}
\newcommand{\lndd}{\Lambda_{n,d}^{*cop}}
\newcommand{\ot}{\otimes}
\renewcommand{\land}{\Lambda_{n,d}}
\newcommand{\dland}{\mathcal{D}(\Lambda_{n,d})}
\newcommand{\set}[1]{\{ #1 \}}
\newcommand{\rep}[1]{\langle{#1}\rangle}
\newcommand{\core}{\operatorname{core}\nolimits}

\newcommand{\evmp}[3]{\ensuremath{M_{2 #1}^+(#2,#3)}}
\newcommand{\evmpm}[3]{\ensuremath{M_{2 #1}^-(#2,#3)}}
\newcommand{\band}[4]{\ensuremath{C_{#1}^{#2}(#3,#4)}}

\makeatletter

\title[Stable Green ring of Taft algebras]{Stable Green ring of the Drinfeld doubles of the generalised Taft algebras (corrections and new results)}
\author[K. Erdmann]{Karin Erdmann}
\address{K. Erdmann, Mathematical Institute,
University of Oxford,
Andrew Wiles Building,
ROQ,
Woodstock Road,
Oxford
OX2 6GG, United Kingdom}
\email{erdmann@maths.ox.ac.uk}
\author[E.L. Green]{Edward L. Green}
\address{E.L. Green, Mathematics Department, Virginia Tech,
Blacksburg, VA  24061-0123, USA}
\email{resolme@math.vt.edu}
\author[N. Snashall]{Nicole Snashall}
\address{N. Snashall, Department of Mathematics, University of Leicester, University Road, Leicester LE1 7RH, United Kingdom}
\email{njs5@le.ac.uk}
\author[R. Taillefer]{Rachel Taillefer}
\address{R. Taillefer, Universit\'e Clermont Auvergne, CNRS, LMBP,  F-63000 Clermont-Ferrand, France}
\email{Rachel.Taillefer@uca.fr}

\date{\today}
\begin{document}

\subjclass[2010]{
17B37, 
81R50, 
16E30, 
06B15, 
16T05. 
}

\keywords{Hopf algebra, Drinfeld double, fusion rules, tensor product, indecomposable module, stable Green ring}

\maketitle

\begin{abstract} \normalsize We return to the fusion rules for the Drinfeld double of the duals of the generalised Taft algebras that we studied in \cite{EGST}. We first correct some proofs and statements in \cite{EGST} that were incorrect, using stable homomorphisms. We then complete this with new results on fusion rules for the modules we had not  studied in \cite{EGST} and a classification of endotrivial and algebraic modules.
\end{abstract}
\section*{Introduction}\label{section:introduction}

Fusion rules, that is, the decomposition of tensor products of modules as a direct sum of indecomposable modules over a Hopf algebra, have been studied in several contexts, such as quantum groups or conformal field theory. In the case of quasi-triangular Hopf algebras, the tensor product of modules over the base field is commutative. Examples of quasi-triangular Hopf algebras are given by Drinfeld doubles of finite dimensional algebras. We are interested here in the Drinfeld doubles $\dland$ of the duals $\Lambda_{n,d}$  of the extended Taft algebras over an arbitrary field $k$ whose characteristic does not divide $d$, where $n$ and $d$ are positive integers with $n$ a multiple of $d$. Such quantum doubles were originally defined  by Drinfeld in order to provide solutions to the quantum
Yang-Baxter equation arising from statistical mechanics. The algebra  $\dland$ has the advantage of being relatively small, and the  Hopf subalgebra $\Lambda_{n,d}$  is finite-dimensional
and basic, but the representations of $\dland$  share
properties with finite-dimensional representations of $U(sl_2)$ and variations.
As an algebra, $\dland $ is tame, and a parametrisation of the
indecomposable modules is known.
Therefore we studied  direct sum decompositions for tensor
products of indecomposable modules in
\cite{EGST}. This has also been done for the Drinfeld double of the Taft algebras, that is, the case $n=d$, in \cite{CMS}, using different methods.

\sloppy

We returned to this problem because  H.-X. Chen
(one of the authors of \cite{CMS}) asked about the proof of \cite[Theorem 4.18]{EGST} (which is wrong) and the proof and statement of \cite[Theorem 4.22]{EGST} (they are both wrong). We are grateful to  H.-X. Chen for drawing our attention to these problems. Moreover, in looking at the details, we also noticed that \cite[Proposition 3.2]{EGST} is incorrect, the proof of \cite[Proposition 4.17]{EGST} is not quite complete, and  there are redundancies in our classification of $\dland$-modules. We also realised that many of the proofs can be simplified by working over a specific block of $\dland$.

In this paper we present a new and more homological approach to these tensor product calculations,
which is based on exploiting stable module homomorphisms. This enables us to provide corrections to the proofs and statements mentioned above (except \cite[Proposition 3.2]{EGST}, which was just a tool we used for some of our  results in \cite{EGST}). Using this new method, we are also able to give general formulas for the decompositions of tensor products involving the remaining modules of even length, the band modules, for which we had only given some examples in \cite{EGST}. As a consequence, we now have a complete description of the stable Green ring of $\dland$, which we give in Section \ref{sec:summary} (Table \ref{table:summary} on page \pageref{table:summary}).

We also include new results, classifying endotrivial and algebraic $\dland $-modules.  If $H$ is a finite-dimensional ribbon Hopf algebra (see for instance \cite[Section 4.2.C]{CP}), then a finite-dimensional $H$-module $M$ is endotrivial if there is an isomorphism $M\otimes_k M^*\cong k \oplus P$ where
$k$ is the trivial $H$-module and $P$ is projective.
Tensoring with an endotrivial module induces an equivalence of the stable module
category, and such equivalences  form  a subgroup of the auto-equivalences of the stable module category. When $H=kG$ is the group algebra of a finite group $G$, endotrivial modules have been studied extensively. They have also been studied for finite group schemes in \cite{CN1,CN2}. However, we have not seen any results on endotrivial modules for other Hopf algebras.
When the Hopf algebra $H$ is our Drinfeld double $\dland$, we  show
that the indecomposable endotrivial modules are precisely the syzygies of the simple modules of dimension $1$ or $d-1$, see Proposition \ref{prop:endotrivial modules}.

The concept of an algebraic module is quite natural; it was introduced as a $kG$-module satisfying a polynomial equation with coefficients in $\zz$ in the Green ring of $kG$, for a finite group $G$. We shall use the following equivalent definition:   a $kG$-module $M$ is called algebraic if the number of non-isomorphic indecomposable summands  of the set of modules $M^{\otimes_k t}$, when $t\pgq 1$ varies, is finite. Such modules occur in particular in the study of the Auslander-Reiten quiver of $kG$. For a study of algebraic $kG$-modules, see for instance \cite{C} and the references there. Here,  we replace the Hopf algebra $kG$ with $\dland$, and we classify the indecomposable algebraic $\dland$-modules.

The paper is organised as follows. In Section \ref{section:list modules}, we describe the quiver and the representations of $\dland$, in particular removing the redundancies mentioned above, and we recall the decomposition of the tensor product of two simple modules. Section \ref{section:proof 4.18} contains our new proof of \cite[Theorem 4.18]{EGST} in Theorem \ref{thm:even tensor simple}. Proposition \ref{prop:case0} is a special case of \cite[Proposition 4.17]{EGST} (enough for our purposes) whose proof is now complete.
Section \ref{section:proof 4.22}  contains a corrected statement and our new proof of \cite[Theorem
4.22]{EGST}, this is  Theorem \ref{thm:even tensor even}. The proof is by induction, the initial
step being Proposition \ref{prop:reduction}. In Section \ref{section:band modules}, we use these
methods to determine the decomposition of the tensor product of a band module with any other
indecomposable module in Propositions \ref{prop:band tensor simple} and \ref{prop:tensor different
  bands} and in Theorem \ref{thm:band tensor band}. We conclude this section with a result that
gives a more conceptual view of the approach we have used above. 
In Section \ref{sec:summary}, we combine all our results on tensor products of representations of
$\dland$ in order to describe the stable Green ring of $\dland$.
Finally, in Section \ref{section:endotrivial}, we present the classification of endotrivial and algebraic $\dland$-modules.

Throughout the paper, $k$ is an algebraically closed field and $q$ is a fixed primitive $d$-th root of unity in $k$ (in particular, the characteristic of $k$ does not divide $d$).  All modules are left modules. The cyclic group of order $n$ is denoted by $\zz_n$ and $\otimes$ denotes the tensor product over $k$.

We refer to \cite{EGST} for all the background on the representation theory of $\dland$, where $n$ and $d$ are integers such that $d$ divides $n$. However, we recall (mainly in Section \ref{section:list modules}) the definitions and results from \cite{EGST} that are required for the understanding of this paper, so that the reader does not need to refer to our original paper.

\section{Parametrization of modules over $\dland$}\label{section:list modules}

In this section, we recall briefly a description of the algebras $\lnd$, $\lndd$ and $\dland$, as well as  the isomorphism classes of representations of $\dland$. We refer to \cite[Sections 2 and 3]{EGST} for more details.

\bn As in \cite[Notation 2.5 and Definition 2.11]{EGST}, we denote by $\langle r \rangle$  the residue of an integer $r$ modulo $d$
taken in the set $\{ 1, 2, \ldots, d\}$ and, for any $u\in\zz_n$, we define a permutation $\sigma_u$ of $\zz_n$ by
\[\sigma_u(j) = d+j-\langle 2j+u-1\rangle\]
 (recall that $d$ divides $n$). If $2j+u-1$ is not divisible by $d$ then the orbit of $j$ under $\sigma_u$
has size $2\frac{n}{d}$ and moreover we have $\sigma_u^{2t}(i)=i+td$ and $\sigma_u(i)^{2t+1}(i)=\sigma_u(i)+td$ in $\zz_n$.
\en

\subsection{The original algebras}\label{subsec:algebras}

The algebra $\lnd$ is described by quiver and relations. Its quiver is the cyclic quiver with vertices $e_1,\ldots,e_{n-1}$  and arrows $a_0,\ldots,a_{n-1}$, where the indices are viewed in $\zz_n$ and each arrow $a_i$ goes from $e_i$ to $e_{i+1}.$ The ideal of relations of $\lnd$ is generated by the paths of length $d$. We also denote by $\gamma_i^m=a_{i+m-1}\cdots a_{i+1}a_i$ the path of length $m$ that starts at $e_i$.

Since we assume that $d$ divides $n$, this algebra is a Hopf algebra by \cite{CC}. We shall only need the antipode here, which is determined by 
\[ S(e_i)=e_{-i}\text{ and }S(a_i)=-q^{i+1}a_{-i-1}\text{ for all } i\in\zz_n. \]

The co-opposite of the dual Hopf algebra, $\lndd$, is the extended Taft algebra, and it is presented by generators and relations:
\[ \lndd=\rep{G,X\,|\, G^n=1,\ X^d=0,\ GX=q^{-1}XG} .\] Its antipode is determined by $S(G)=G^{-1}$ and $S(X)=-XG^{-1}$.

These algebras are both Hopf subalgebras of the Drinfeld double $\dland$, that is equal to $\land\ot\lndd$ as a vector space, and has a basis given by the set of $G^iX^j\gamma_\ell^m$ with $i,\ell$ in $\zz_n$ and $0\ppq j,m\ppq d-1$.  In this paper we shall not need the relations between the generators.

\subsection{$\dland$-modules of odd length, projective $\dland$-modules and blocks of $\dland$}\label{subsec:odd}

The simple $\dland$-modules are labelled  $L(u,i)$ for $(u, i) \in \zz_n^2$. The description below is taken from \cite[Section 2, mainly Propositions 2.17 and 2.21]{EGST}.

The module $L(u,i)$ is
projective if and only if $2i+u-1$ is divisible by $d$, in which case $\dim L(u,i)=d$. When $L(u,i)$ is not projective, then $\dim L(u,i)=d-\langle 2i+u-1\rangle=\sigma_u(i)-i$ and $L(u,i)$ contains two distinguished vectors  $\tilde{H}_{u,i}$ and $\tilde{F}_{u,i}$, with
the following properties:
\begin{enumerate}[(a)]
\item $\tilde{H}_{u,i}$ spans the kernel of the action of $X$ on
$L(u,i)$ as a vector space
and $e_j\tilde{H}_{u,i} = \delta_{ij}\tilde{H}_{u,i}$,
\item   $e_j\tilde{F}_{u,i} = \delta_{j,\sigma_u(i)-1}\tilde{F}_{u,i}$ and the element $\tilde{F}_{u,i}$
is annihilated by  all the arrows in the quiver of $\Lambda_{n,d}$,
\end{enumerate} where $\delta$ is the Kronecker symbol.
The module $L(u,i)$ has basis $\{ X^t\tilde{F}_{u,i}: 0\ppq t< \dim L(u,i)=N\}$ and $X^{N-1}\tilde{F}_{u,i}$ is a non-zero scalar multiple of $\tilde{H}_{u,i}$. Moreover, the action of $G$ on these basis elements is given by $GX^{t}\tilde{F}_{u,i}=q^{-i-t}\tilde{F}_{u,i}.$

In \cite[Proposition 2.21]{EGST}, we characterised the simple modules as follows.

\bp \label{prop:simples} Let $S$ be a simple module. Set $E_u=\frac{1}{n}\sum_{i,j\in\zz_n}q^{-i(u+j)}G^ie_j$ for $u\in\zz_n$. Then $S$ is isomorphic
to $L(u,i)$ if and only if the three following properties hold: 

\begin{enumerate}[(a)]
\item  $\dim S=\dim L(u,i).$
\item $E_u$ acts as identity on
  $S$, and $E_v$ acts as zero on $S$ if $v\neq u$. 
\item Let $Y$ be the generator of $S$ which is in the kernel of the
  action of $X$ (this is well-defined up to a non-zero scalar and corresponds
  to $\tilde{H}_{u,i}$). Then the vertex $e_i$ acts as identity on
  $Y$, and the other vertices act as zero. 
\end{enumerate}
 \ep

\bigskip

The projective cover $P(u,i)$ of $L(u,i)$ has four composition factors. The socle and the top are isomorphic to $L(u,i)$, and
\[\rad P(u,i)/\soc P(u,i) \cong L(u, \sigma_u(i)) \oplus L(u, \sigma_u^{-1}(i)).\]

\bigskip

The \textbf{indecomposable modules of odd length} are precisely the syzygies of the non-projective simple modules, that is, the   $\Omega^m(L(u,i))$ for $m\in\zz$ and $(u,i)\in\zz_n^2$ where $d$ does not divide $2i+u-1$.

\bigskip

It follows from the structure of the projective modules that the simple modules in the block of $L(u,i)$ are precisely the simple
modules $L(u, \sigma_u^t(i))$  for all $t$, and there are
$2\frac{n}{d}$ of them. Theorem 2.26 of \cite{EGST} completely describes the basic algebra
of a non-simple block.
Each block is symmetric and special biserial with radical cube zero.

More precisely, as an algebra, $\dland$ is the direct sum of $\frac{n^2}{d}$ simple blocks and of $\frac{n(d-1)}{2}$ blocks $\Bb_{u,i}$ for $(u,i)\in\zz_n^2$ such that $2i+u-1$ is not a multiple of $d$. The quiver of $\Bb_{u,i}$ is 
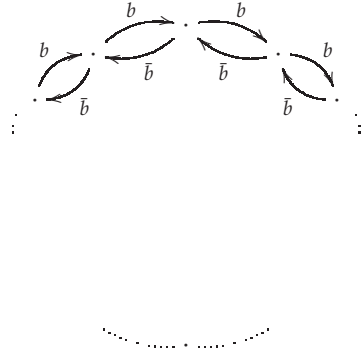
\[\xymatrix@=.01cm{
&&&&&&&&&&&\cdot\ar@/^.5pc/[rrrrrd]^{b}\ar@/^.5pc/[llllld]^{\bar{b}}\\
&&&&&&\cdot\ar@/^.5pc/[rrrrru]^{b}\ar@/^.5pc/[llldd]^{\bar{b}}&&&&&&&&&&\cdot\ar@/^.5pc/[rrrdd]^{b}\ar@/^.5pc/[lllllu]^{\bar{b}}\\\\
&&&\cdot\ar@/^.5pc/[rrruu]^{b}\ar@{.}@/_.3pc/[ldd]&&&&&&&&&&&&&&&&\cdot\ar@/^.5pc/[llluu]^{\bar{b}}\ar@{.}@/^.3pc/[rdd]\\\\
&&&&&&&&&&&&&&&&&&&&&&&\\
\\
\\
\\\\\\\\\\\\\\\\
&&&&&&&&&&&&&&&&\\
&&&&&&&&&&&\cdot\ar@/_.3pc/@{.}[rrrrru]\ar@/^.3pc/@{.}[lllllu]
}\] with $\frac{2n}{d}$ vertices and $\frac{4n}{d}$ arrows. The
relations on this quiver are $bb$, $\bar{b}\bar{b}$ and
$b\bar{b}-\bar{b}b$ (there are $\frac{6n}{d}$ relations on each of
these quivers). The vertices in this quiver correspond to the simple modules $L(u,i),$
$L(u,\sigma_u(i)),$ $L(u,\sigma_u^2(i)),$ $\ldots,$ $L(u,\sigma_u^{\frac{2n}{d}-1}(i)).$ Hence $\Bb_{u,i}=\Bb_{v,j}$ if and only if $u=v$ and $j=\sigma_u^t(i)$ for some $t\in\zz.$

\medskip

Moreover, the arrows $b$ and $\bar{b}$ are described as follows. For each $p$, there are basis elements in $L(u,\sigma^p(i))$ that, following the notation in \cite[Propostion 2.17]{EGST}, we denote by $\tilde{D}_{u,\sigma_u^{p+1}(i)}$ and $F_{u,\sigma_u^{p-1}(i)}$, such that:
\begin{itemize}
\item the action of $b_p$ on $L(u,\sigma^p(i))$ is given by multiplication by $\gamma_{\sigma_u^p(i)}^{\dim L(u,\sigma^{p}_u(i))}$ and $\gamma_{\sigma_u^{p}(i)}\tilde{D}_{u,\sigma_u^{p+1}(i)}$ is a non-zero scalar multiple of $\tilde{H}_{u,\sigma_u^{p+1}(i)}$,
\item the action of $\bar b_{p-1}$ on $L(u,\sigma^p(i))$ is given by multiplication by $X^{\dim L(u,\sigma^p_u(i))}$ and $X^{\dim L(u,\sigma^p_u(i))}F_{u,\sigma_u^{p-1}(i)}$ is equal to $\tilde{F}_{u,\sigma_u^{p-1}(i)}$.
\end{itemize}

\subsection{$\dland$-modules of even length}\label{subsec:even}

We described all the non-projective indecomposable representations of $\dland$ in \cite[Section 3 and Appendix A]{EGST}. These are
\begin{enumerate}[(a)]
\item The modules of odd length described above.
\item The string modules of  length $2\ell$, denoted by $M_{2\ell}^+(u,i)$ and $M_{2\ell}^-(u,i)$, for any positive integer $\ell$ and  $(u,i)\in  \zz_n^2$ with $2i+u-1\not\equiv0\pmod d$.
\item The band modules $C_\lambda^{\ell}(u,i)$ of length $2\ell\frac{n}{d}$ for $\lambda\in k\setminus\left\{0\right\}$, any positive integer $\ell$ and $(u,i)\in  \zz_n^2$ with $2i+u-1\not\equiv0\pmod d$, taking $i$ modulo  $d$ (that is, one
for each orbit of the square of $\sigma_u$).
\end{enumerate}

We now describe the modules of even length in more detail for future use.
Fix a block $\Bb_{u,i}.$

\subsubsection{String modules of even length.}
 For each $0\ppq p\ppq \frac{2n}{d}-1$ and for each $\ell \pgq 1,$ there are two indecomposable
  modules of length $2\ell$ which we call $M^{\pm }_{2\ell}(u,\sigma_u^p(i)):$ 
\begin{enumerate}[$\bullet$]
\item The module $M^{+}_{2\ell}(u,i)$ has top composition factors
  $L(u,i),$
  $L(u,\sigma_u^{2}(i)),$ $\ldots,$ $L(u,\sigma_u^{2(\ell-1)}(i))$ and
    socle composition factors $L(u,\sigma_u(i)),$ $
  L(u,\sigma_u^{3}(i)),$ $\ldots,$
  $L(u,\sigma_u^{2(\ell-1)+1}(i))$:
\[\def\objectstyle{\scriptstyle}\xymatrix@=.05cm{L(u,i) \ar@{-}[dr]  &  &
    L(u,\sigma_u^{2}(i))\ar@{-}[dr]\ar@{-}[dl] &  & 
  L(u,\sigma_u^{2(\ell-1)}(i))\ar@{-}[dr]\ar@{-}[dl]&&              
 \\&L(u,\sigma_u(i)) && \cdots &&
L(u,\sigma_u^{2(\ell-1)+1}(i)) &&     
}\] The lines joining the simple modules are given by multiplication
by the appropriate $b$-arrow or $\bar{b}$-arrow (in the case $n=d$,
when there is an ambiguity, the first line is multiplication by
$\gamma^{\dim L(u,i)}$, the next one is multiplication by
a non-zero scalar multiple of $X^{d-\dim L(u,i)}$, and so on, up to non-zero scalars).

\item The module $M^{-}_{2\ell}(u,i)$ has top composition factors
  $L(u,i),$ $
  L(u,\sigma_u^{-2}(i)),$ $\ldots,$ $L(u,\sigma_u^{-2(\ell-1)}(i))$ and
    socle composition factors $L(u,\sigma_u^{-1}(i)),$ $
  L(u,\sigma_u^{-3}(i)),$ $\ldots,$
  $L(u,\sigma_u^{-2(\ell-1)+1}(i))$:
$$\def\objectstyle{\scriptstyle}\xymatrix@=.05cm{ & L(u,\sigma_u^{-2(\ell-1)}(i))\ar@{-}[dr]\ar@{-}[dl]  &&  L(u,\sigma_u^{-2}(i)) \ar@{-}[dr] \ar@{-}[dl] &&
L(u,i) \ar@{-}[dl]  
\\L(u,\sigma_u^{-2(\ell-1)+1}(i)) &&   \ldots && L(u,\sigma_u^{-1}(i))
} $$ As for the other string modules, the lines represent
multiplication by an appropriate $b$ or $\bar{b}$ arrow, and in the
case $n=d$ the first one from the left is multiplication by a power of $X$ and so on.
\end{enumerate} In both cases, indices are taken modulo $\frac{2n}{d}.$

These modules are periodic of period $\frac{2n}{d}.$ Moreover, the Auslander-Reiten sequences of the string modules $M_{2\ell}^{\pm}(u,i)$ are given by
\begin{align*}
&\A(M_{2\ell}^+(u,i)):\qquad0\rightarrow M_{2\ell}^+(u,i-d)\rightarrow M_{2\ell+2}^+(u,i-d)\oplus  M_{2\ell-2}^+(u,i) \rightarrow M_{2\ell}^+(u,i)\rightarrow 0,\\
&\A(M_{2\ell}^-(u,i)):\qquad0\rightarrow M_{2\ell}^-(u,i+d)\rightarrow M_{2\ell+2}^-(u,i+d)\oplus  M_{2\ell-2}^-(u,i) \rightarrow M_{2\ell}^-(u,i)\rightarrow 0
\end{align*} (where $M_0^{\pm}(u,i)=0$).

\subsubsection{Band modules.} For
  each $\lambda\neq 0$ in $k,$ and for each $\ell\pgq 1$, there is an
  indecomposable module of length $\frac{2n}{d}\ell$, which we denote
  by $C^{\ell}_\lambda(u,i)$. It is
  defined as follows. Denote by $\epsilon_p$, for $p\in\zz_{2n/d}$, the vertices in the quiver of $\Bb_{u,i}$, with $\epsilon_0$ corresponding to $L(u,i)$. The arrow $b_p$ goes from $\epsilon_p$ to $\epsilon_{p+1}$ and the arrow $\bar{b}_{p}$ goes from $\epsilon_{p+1}$ to $\epsilon_{p}$.  Let $V$ be an $\ell$-dimensional vector space. Then
  $C^{\ell}_\lambda(u,i)$ has underlying space
  $C=\bigoplus_{p=0}^{\frac{2n}{d}-1}C_p$ with $C_p=V$ for all $p.$ The
  action of the idempotents $\epsilon_p$ is such that
  $\epsilon_pC=C_p.$ The action of the arrows $\bar{b}_{2p}$ and
  $b_{2p+1}$ is zero. The action of the arrows $\bar{b}_{2p+1}$ is the
  identity of $V.$ The action of the arrows $b_{2p}$ with $p\neq 0$ is
  also the identity. Finally, the action of $b_0$ is given by the
  indecomposable Jordan matrix $J_\ell(\lambda).$

It is periodic of period 2 and its Auslander-Reiten sequence is given by 
\[ \A(C_\lambda^{\ell}(u,i)):\qquad0\rightarrow C_\lambda^{\ell}(u,i-d)\rightarrow C_\lambda^{\ell+1}(u,i)\oplus  C_\lambda^{\ell-1}(u,i) \rightarrow C_\lambda^{\ell}(u,i)\rightarrow 0 \] (where $C_\lambda^0(u,i)=0$).

Note that $\soc(C)=\rad(C)=\bigoplus_{p}{\epsilon_{2p+1}C}$ and that
$C/\rad(C)=\bigoplus_{p}{\epsilon_{2p}C}$.

\br
In \cite{EGST}, we denoted  $C^{\ell}_\lambda(u,i)$ by  $C^{\ell+}_\lambda(u,i)$, and we had more band modules. First note that  $C_{\lambda}^{\ell}(u,i)$ is isomorphic to
$C_{\lambda}^{\ell}(u, \sigma_u^2(i))$, so that in (c) above we have removed some repetitions that occurred in \cite{EGST}.

In addition, in \cite{EGST} we also had the modules $C_\lambda^{\ell -}(u,i)$, defined in a similar way by interchanging  $b$'s and $\bar{b}$'s.
However, we should have noted that for each $\lambda\in k\setminus\left\{0\right\}$, there exists $\mu\in k\setminus\left\{0\right\}$ and $(v,j)\in\zz_n^2$ such that $C_\lambda^{\ell -}(u,i)\cong C_\mu^{\ell}(v,j)$.

Indeed, by definition of the module $C^{\ell -}_\lambda(u,i)$, the action of $\bar{b}_0$ is given by $J_\ell(\lambda)$, the action of $\bar{b}_i$ for $i$ even, $i\neq0$, is given by $\id$, the action of $b_i$ for $i$ odd is given by $\id$ and the action of the other arrows is given by $0.$ Changing bases of the vector spaces $C_p^-$, we can ensure that the action of $b_{-1}$ is given by $J_\ell\left(\frac{1}{\lambda}\right) $, the action of $\bar{b}_0$ is given by $\id$ and the rest is unchanged, thus obtaining the module $C^{\ell}_{\frac{1}{\lambda}}(u,\sigma_u^{-1}(i))$. Hence $C^{\ell -}_\lambda(u,i)\cong C^{\ell}_{\frac{1}{\lambda}}(u,\sigma_u^{-1}(i))$.

Therefore we need only consider one family of band modules, which we denote by $C_\lambda^{\ell}(u,i)$ for $\lambda\in  k\setminus\left\{0\right\}$, and  we need only one of these for each orbit of $\sigma_u^2$ in $\zz_n$.

However, it should be noted that this parameter $\mu$ such that  $C_\lambda^{\ell -}(u,i)\cong C_\mu^{\ell}(v,j)$ is not well defined, unless we have fixed the block $\Bb_{u,i}$ (including which vertex is labelled $\epsilon_0$): using the definition in \cite[Section A.2.2]{EGST}, where the representations of the basic algebra of a given block are described, the module $C_\lambda^{\ell-}(u,i)$ is  defined by considering that the block $\Bb_{u,i}$ is equal to $\Bb_{u,\sigma_u(i)}$ with a shift in the indices of the vertices and arrows by $1$; here we get immediately that $C_\lambda^{\ell -}(u,i)=C_\lambda^{\ell }(u,i)$.

\er

As we mentioned above, there is a difficulty in differentiating modules of length $2\frac{n}{d}$, since $M_{2\frac{n}{d}}^+(u,i)$, $M_{2\frac{n}{d}}^-(u,i)$ and $C_\lambda^{1}(u,i)$ are modules of length $2\frac{n}{d}$ with the same top and the same socle. The following property allows us to distinguish them. Recall that $\lnd$ and $\lndd$ are subalgebras of $\dland$, so that any $\dland$-module can be viewed as a $\lnd$-module and as a $\lndd$-module.

\begin{pty}\label{pty:charact band length pullback}
Assume that $d$ does not divide $2i+u-1.$ Then
\begin{itemize}
\item the modules $M_{2\ell}^+(u,i)$ are projective as $\lnd$-modules but not as $\lndd$-modules,

\item the modules $M_{2\ell}^-(u,i)$ are projective as $\lndd$-modules but not as $\lnd$-modules,

\item the modules $C_\lambda^\ell(u,i)$ are projective both as $\lnd$-modules and as $\lndd$-modules.
\end{itemize}
\end{pty}

\bpf We use the notation from the last part of Subsection \ref{subsec:odd}.

It follows from \cite{CC} and \cite{LZ} that the indecomposable projective modules over $\lnd$ or $\lndd$
 are precisely the indecomposable modules of dimension $d$.
Therefore $M_2^+(u,i)$, which is equal to $\lnd \tilde{D}_{u,\sigma_u(i)}$ as a $\lnd$-module,  is indecomposable of dimension $d$, hence projective. Moreover, as a $\lnd$-module,
 $M_{2\ell}^+(u,i)\cong\bigoplus_{t=0}^{d-1} M_2^+(u,i+td)$ is also projective.

Similarly, the $\lndd$-module $M_2^-(u,i)=\lndd F_{u,\sigma_u^{-1}(j)}$ is projective, and $M_{2\ell}^-(u,i)\cong \bigoplus_{t=0}^{\ell-1} M_2^-(u,i-td)$ is a projective  $\lndd$-module.

 As a  $\lndd$-module, $M_{2\ell}^+(u,i)$  decomposes as $ L(u,i)\oplus\bigoplus_{t=1}^{\ell-1}M_2^-(u,i+td)\oplus L(u,\sigma_u(i)+(\ell-1)d)$, therefore it has a  summand $L(u,i)$ whose dimension is strictly less than $d$, hence that is not projective. Therefore $M_{2\ell}^+(u,i)$ is not projective.

Similarly, $M_{2\ell}^-(u,i)$ is not projective as a $\lnd$-module.

Finally, as a $\lnd$-module, $C_\lambda^\ell(u,i)\cong \bigoplus_{t=0}^{\frac{n}{d}-1} M_2^+(u,i+td)^\ell$ is projective, and
as a $\lndd$-module,  $C_\lambda^\ell(u,i)\cong \bigoplus_{t=0}^{\frac{n}{d}-1} M_2^-(u,i-td)^\ell$ is projective.
\epf

\br\label{rk:tensor restricted projective} It is clear that the $\dland$-modules $P(u,i)$ are projective both as $\lnd$-modules and as $\lndd$-modules. Moreover, the non-projective simple $\dland$-modules $L(u,i)$ have dimension at most $d-1$ and therefore they are  not projective as $\lnd$-modules or as $\lndd$-modules.

In this paper, we consider tensor products of $\dland$-modules over the base field. The tensor products of $\lnd$-modules and of $\lndd$-modules have been studied in  \cite{CC} and \cite{LZ}. It follows in particular from their results that  if $M$ and $N$ are two indecomposable $\lnd$-modules (respectively $\lndd$-modules), then $M\otimes N$ is projective if and only if $M$ or $N$ is projective.
\er

\bd We shall say that two string modules $M_{2\ell}^{+}(u,i)$ and $M_{2t}^+(u,j)$ (respectively  two string modules  $M_{2\ell}^{-}(u,i)$ and $M_{2t}^-(u,j)$, respectively two band modules $C_\lambda^\ell(u,i)$ and $C_\mu^t(u,j)$) have the \emph{same parity} if  $i\equiv j\pmod{d}$ and that they have \emph{different parities} if $j\equiv \sigma_u(i)\pmod d$.
\ed

This definition is consistent with \cite[Definition A.3]{EGST}.

\bd We shall  write $\core(M)$ for the direct sum of the non-projective indecomposable summands in $M$.
\ed

In  \cite[Theorem 4.1]{EGST}, we determined the  decomposition of the tensor product of two simple modules. Here we shall need the core of this module.

\bp(cf. \cite[Theorem 4.1]{EGST})\label{prop:simple tensor simple} Let $L(u,i)$ and $L(v,j)$ be two non-projective simple $\dland$-modules. Then
\[ \core(L(u,i)\otimes L(v,j))\cong\bigoplus_{\theta\in \I}L(u+v,i+j+\theta) \] where $\I=
\begin{cases}
\set{\theta\,;\, 0\ppq \theta\ppq \min(\dim L(u,i),\dim L(v,j))-1}\text{ if }\dim L(u,i)+\dim L(v,j)\ppq d\\
\set{\theta\,;\, \dim L(u,i)+\dim L(v,j)-d\ppq \theta\ppq \min(\dim L(u,i),\dim L(v,j))-1}\text{ otherwise. }
\end{cases}
$

Moreover, if $\dim L(u,i)+\dim L(v,j)\ppq d$, we have $L(u,i)\otimes L(v,j)=\core(L(u,i)\otimes L(v,j)).$
\ep

\bpf Set $N=\dim L(u,i)$ and $N'=\dim L(v,j)$. Without loss of generality, we can assume that $N\ppq N'.$ We proved in \cite[Theorem 4.1]{EGST} that
\[ L(u,i)\otimes L(v,j)\cong \begin{cases} \bigoplus_{\theta=0}^{N-1}L(u+v,i+j+\theta)&\text{ if }N+N'\ppq d+1\\
 \bigoplus_{\theta=N+N'-d}^{N-1}L(u+v,i+j+\theta)\oplus \text{ projective }&\text{ if }N+N'\pgq d+1
\end{cases} \] so we need only determine which of the simple modules are projective. Note that $2i+u-1\equiv -N\pmod d$ and $2j+v-1\equiv -N'\pmod d$, therefore $2(i+j+\theta)+(u+v)-1\equiv -N-N'+2\theta+1\pmod d.$ Moreover, since $\theta\ppq N-1\ppq N'-1$, we have $ -N-N'+2\theta+1\ppq -1.$

If $N+N'\ppq d$ we have $ -N-N'+2\theta+1\pgq -d+1$, and if $N+N'>d+1$ we have $-N-N'+2\theta+1\pgq -N-N'+2(N+N'-d)+1=N+N'-2d+1>d+1-2d+1=2-d$, therefore in both cases $2(i+j+\theta)+(u+v)-1\not\equiv 0\pmod d$ and $L(u+v,i+j+\theta)$ is not projective.

If $N+N'=d+1$, we have $-N-N'+2\theta+1\pgq 2\theta-d\equiv 2\theta\pmod d$ with $0\ppq 2\theta\ppq 2(N-1)\ppq N+N'-2=d-1$, therefore $2(i+j+\theta)+(u+v)-1\equiv 0\pmod d$ if and only if $\theta=0$.

Finally, the only projective that can occur is $L(u+v,i+j)$ when $N+N'=d+1$. The result then follows.
\epf

\section{Decomposition of the tensor product of a simple module with a string module of even length} \label{section:proof 4.18}

The stable Green ring $r_{st}(H)$ of a finite dimensional Hopf algebra $H$ is the ring whose $\zz$-basis is the set of isomorphism classes of indecomposable modules in the stable category $H$-$\smod$,
that is, the non-projective indecomposable $H$-modules, and whose addition and multiplication are induced respectively by the direct sum and the tensor product.

The aim of this section is to determine the tensor product of a string module of even length with a simple module up to projectives, from which we can deduce  the tensor product of a string module of even length with a string module of odd length in the stable  Green ring $r_{st}(\dland)$. This will provide a correction to the proof of \cite[Theorem 4.18]{EGST}.

For the proof,  we will use a general result involving splitting trace
modules.
We  recall that a module $M$ is a \emph{splitting trace module} if the trivial module $L(0,0)$ is a direct summand in $\End_k(M)$. Moreover, for any ribbon Hopf algebra such as $\dland$, there are isomorphisms $\End_k(M)\cong M^*\ot M\cong M\ot M^*$ of $\dland$-modules for any module $M$ (where $M^*$ is the $k$-dual of $M$), see for instance \cite[Section 4.2.C]{CP}.

\medskip

The \textbf{indecomposable splitting trace modules for $\dland$} were given in \cite[Subsection 3.3]{EGST}: they are precisely the non-projective indecomposable modules of odd length.

\bl\label{lemma:AC} Let $C$ be an indecomposable non-projective
$\dland$-module and let $\A(C)$ be its Auslander-Reiten sequence. For any $\dland$-module $M$, the following are equivalent.
\begin{enumerate}[(i)]
\item $C$ is a direct summand in $\End_k(M)\ot C$.
\item $C$ is a direct summand in $M^*\otimes M\ot C$.
\item The sequence $\A(C)\ot M$ does not split.
\item If  $0\rightarrow A\rightarrow B\rightarrow C\rightarrow0$ is any non-split exact sequence of $\dland$-modules, then the exact sequence $0\rightarrow A\otimes M\rightarrow B\otimes M\rightarrow C\otimes M\rightarrow0$ does not split.
\end{enumerate} If moreover $M$ is a splitting trace module, then all these properties hold.
\el

\bpf The equivalence of the first two statements follows from the isomorphism $\End_k(M)\cong M^*\ot M$ of $\dland$-modules. The equivalence between statements (ii), (iii) and (iv) is proved in the same way as  \cite[Proposition 2.3]{AC}. Finally, if $M$ is a splitting trace module, then $L(0,0)$ is a summand in $\End_k(M)$ so that $C\cong L(0,0)\otimes C$ is a direct summand in $\End_k(M)\otimes C$ and (i) holds.
\epf

We start with the special cases of the tensor product of $M_2^{\pm}(0,td)$ with a simple module in Proposition \ref{prop:special M2 tensor L} and of  the tensor product of $M_{2\ell}^{\pm}(0,td)$ with a simple module in Proposition \ref{prop:case0}, from which we deduce the general case in Theorem \ref{thm:even tensor simple}.

\bp (cf. \cite[Proposition 4.17]{EGST}) \label{prop:special M2 tensor L} For all $t\in\zz,$ we have $\core(M_2^+(0,td)\otimes L(v,j))\cong M_2^+(v,j+td)$ and $\core(M_2^-(0,td)\otimes L(v,j))\cong M_2^-(v,j+td)$.
\ep

\bpf   Set $N=\dim L(v,j).$  We tensor the non-split exact sequence $0\rightarrow L(0,1+td)\rightarrow M_2^+(0,td)\rightarrow L(0,td)\rightarrow 0$ by $L(v,j)$. Using   Proposition \ref{prop:simple tensor simple}  as well as  $\dim L(0,td)=1$ and $\dim L(0,1+td)=d-1$ (so that $\#\I=1$), we obtain  an exact sequence $0\rightarrow L(v,j+td+N)\oplus P \rightarrow M_2^+(0,td)\otimes L(v,j)\rightarrow L(v,j+td)\rightarrow 0$ with $P$ a projective module.
Since $P$ is injective, we have $M_2^+(0,td)\otimes L(v,j)\cong U\oplus P$ for some module $U$ and we obtain an exact sequence $0\rightarrow L(v,j+td+N)\rightarrow U\rightarrow   L(v,j+td)\rightarrow 0$.
The sequence cannot split by Lemma \ref{lemma:AC} because $L(u,j)$ is a splitting trace module,
 so $U$ is an indecomposable module of length $2$ with  $\Top U=L(v,j+td)$ and $\soc U=L(v,j+td+N)=L(v,\sigma_v(j+td))$. Moreover, it follows from Property \ref{pty:charact band length pullback} and Remark \ref{rk:tensor restricted projective} that $M_2^+(0,td)\otimes L(v,j)$, and hence $U$, is projective as a $\lnd$-module but not as a $\lndd$-module, therefore  we must have $U=M_{2}^+(v,j+td).$

The proof of the second part is similar.
\epf

In the next lemma as well as in Proposition \ref{prop:reduction}, we shall use properties of homomorphisms between string modules of even length over a given block; these properties can easily be seen by working over the basic algebra of the block.

\bl\label{lemma:dim Ext1} Let $t$ and $\ell$ be integers with $1\ppq t\ppq \ell.$ Then the spaces $\Ext^1_{\dland}(M_{2\ell}^+(u,i+ d),M_{2t}^+(u,i))$ and $\Ext^1_{\dland}(M_{2\ell}^-(u,i- d),M_{2t}^-(u,i))$ have dimension
\[ \#\set{0\ppq y\ppq t-1\,|\, y\equiv 0\pmod{\frac{n}{d}} }-\#\set{0\ppq y\ppq t-1\,|\, y\equiv \ell\pmod{\frac{n}{d}} }+\#\set{0\ppq y\ppq t-1\,|\, y\equiv t-\ell-1\pmod{\frac{n}{d}} }. \] In particular, $\Ext^1_{\dland}(M_{2\ell}^+(u,i+ d),M_{2}^+(u,i))$ and  $\Ext^1_{\dland}(M_{2\ell}^-(u,i- d),M_{2}^-(u,i))$ have dimension $1$. A basis of $\Ext^1_{\dland}(M_{2\ell}^+(u,i+ d),M_{2}^+(u,i))$ (respectively $\Ext^1_{\dland}(M_{2\ell}^-(u,i- d),M_{2}^-(u,i))$) is given by the equivalence class of the exact sequence
\begin{align*}
&0\rightarrow M_{2}^+(u,i)\rightarrow M_{2\ell+2}^+(u,i)\rightarrow M_{2\ell}^+(u,i+d)\rightarrow 0\\
(respectively\quad &0\rightarrow M_{2}^-(u,i)\rightarrow M_{2\ell+2}^-(u,i)\rightarrow M_{2\ell}^-(u,i- d)\rightarrow 0).
\end{align*}
\el

\begin{proof} We prove the first statement, the proof of the second is similar.

 Applying $\Hom_{\dland}(-,M_{2t}^+(u,i))$ to the exact sequence
\[ 0\rightarrow \Omega(M_{2\ell}^+(u,i+d))\cong M_{2\ell}^+(u,\sigma_u(i))\rightarrow Q:=\bigoplus_{x=1}^{\ell}P(u,i+xd)\rightarrow M_{2\ell}^+(u,i+d)\rightarrow 0 \] gives an exact sequence
\begin{align*}
0&\rightarrow \Hom_{\dland}(M_{2\ell}^+(u,i+d),M_{2t}^+(u,i)) \rightarrow \Hom_{\dland}(Q,M_{2t}^+(u,i)) \rightarrow \Hom_{\dland}(M_{2\ell}^+(u,\sigma_u(i)),M_{2t}^+(u,i)) \\&\rightarrow \Ext^1_{\dland}(M_{2\ell}^+(u,i+d),M_{2t}^+(u,i))\rightarrow0
\end{align*}
so that
\begin{align*}
\dim\Ext^1_{\dland}(M_{2\ell}^+(u,i+d),M_{2t}^+(u,i))&=\dim \Hom_{\dland}(M_{2\ell}^+(u,i+d),M_{2t}^+(u,i)) \\&\ -\dim \Hom_{\dland}(Q,M_{2t}^+(u,i)) \\&\ +\dim \Hom_{\dland}(M_{2\ell}^+(u,\sigma_u(i)),M_{2t}^+(u,i)).
\end{align*}

We determine each of these dimensions.

 We have $\dim \Hom_{\dland}(Q,M_{2t}^+(u,i))=\sum_{x=1}^\ell\#\set{y\,;\, 0\ppq y\ppq t-1,\ y\equiv x\pmod{\frac{n}{d}}}.$

 \sloppy Since $M_{2\ell}^+(u,\sigma_u(i))$ and $M_{2t}^+(u,i)$ have different parities, any non-zero map in $\Hom_{\dland}(M_{2\ell}^+(u,\sigma_u(i)),M_{2t}^+(u,i))$ must send the top of $M_{2\ell}^+(u,\sigma_u(i))$ into the socle of $M_{2t}^+(u,i)$. Therefore
\begin{align*}
\dim\Hom_{\dland}(M_{2\ell}^+(u,\sigma_u(i)),M_{2t}^+(u,i))&=\sum_{x=0}^{\ell-1}\#\set{y\,;\, 0\ppq y\ppq t-1,\ y\equiv x\pmod{\frac{n}{d}}}\\&=\dim\Hom_{\dland}(Q,M_{2t}^+(u,i))\\&\qquad +\#\set{y\,;\, 0\ppq y\ppq t-1,\ y\equiv 0\pmod{\frac{n}{d}}}\\&\qquad -\#\set{y\,;\, 0\ppq y\ppq t-1,\ y\equiv \ell\pmod{\frac{n}{d}}}.
\end{align*}

 Since $M_{2\ell}^+(u,i+d)$ and $M_{2t}^+(u,i)$ have the same parity, the image of any non-zero map in $\Hom_{\dland}(M_{2\ell}^+(u,i+d),M_{2t}^+(u,i))$ is a quotient of $M_{2\ell}^+(u,i+d)$ and a string submodule of $M_{2t}^+(u,i)$ that starts at $L(u,i)$. Therefore, since $\ell\pgq t$,
\begin{align*}
\dim\Hom_{\dland}(M_{2\ell}^+(u,i+d),M_{2t}^+(u,i))&=\#\set{x\,;\, \ell-t+1\ppq x\ppq \ell,\ x\equiv 0\pmod{\frac{n}{d}}}\\&=\#\set{y\,;\, 0\ppq y\ppq t-1,\ y\equiv t-\ell-1\pmod{\frac{n}{d}}}.
\end{align*}

The general formula follows from there and the  case $t=1$ is then clear.
\end{proof}

\bp\label{prop:case0} For any integer $\ell\pgq 1$ and all $t\in\zz$, we have
\begin{align*}
\core(M_{2\ell}^+(0,td)\otimes L(v,j))&\cong M_{2\ell}^+(v,j+td)\text{ and }\\
\core(M_{2\ell}^-(0,td)\otimes L(v,j))&\cong M_{2\ell}^-(v,j+td).
\end{align*}
\ep

\begin{proof}
The proof is by induction on $\ell.$ Proposition \ref{prop:special M2 tensor L} shows that the result is true for $\ell=1.$

 Now assume that $\core(M_{2\ell}^+(0,td)\otimes L(v,j))\cong M_{2\ell}^+(v,j+td)$ for all $t\in\zz$ and for a given $\ell\pgq 1.$ There is a non-split exact sequence
\[ 0\rightarrow M_{2}^+(0,td)\rightarrow M_{2\ell+2}^+(0,td)\rightarrow M_{2\ell}^+(0,(t+1)d)\rightarrow0. \] Tensoring with $L(v,j)$ gives an exact sequence
\[ 0\rightarrow M_{2}^+(v,j+td)\oplus P_1\rightarrow M_{2\ell+2}^+(0,td)\otimes L(v,j)\rightarrow M_{2\ell}^+(v,j+(t+1)d) \oplus P_2\rightarrow 0  \] with $P_1$ and $P_2$ projective-injective modules, so that, factoring out the split exact sequence $0\rightarrow P_1\rightarrow P_1\oplus P_2\rightarrow P_2\rightarrow 0$, we have an exact sequence
\[ 0\rightarrow M_{2}^+(v,j+td)\rightarrow \core(M_{2\ell+2}^+(0,td)\otimes L(v,j))\oplus P\rightarrow M_{2\ell}^+(v,j+(t+1)d) \rightarrow 0   \] for some projective module $P.$ Moreover, since $L(v,j)$ is a splitting trace module, by Lemma \ref{lemma:AC} this sequence is not split. By Lemma \ref{lemma:dim Ext1}, we have $P=0$ and $\core(M_{2\ell+2}^+(0,td)\otimes L(v,j))\cong M_{2\ell+2}^+(v,j+td)$ thus proving the induction.

The proof of the isomorphism $\core(M_{2\ell}^-(0,td)\otimes L(v,j))\cong M_{2\ell}^-(v,j+td)$ is similar, starting with the exact sequence \[ 0\rightarrow M_{2}^-(0,td)\rightarrow M_{2\ell+2}^-(0,td)\rightarrow M_{2\ell}^-(0,(t-1)d)\rightarrow0. \qedhere \]
\end{proof}

\bt(cf. \cite[Theorem 4.18]{EGST})\label{thm:even tensor simple} Fix an integer $\ell\pgq 1.$ For any $(u,i)$ and $(v,j)$ in $\zz_n^2$ with $2i+u-1\not\equiv 0\pmod{d}$ and $2j+v-1\not\equiv0\pmod d$, we have
\begin{align*}
\core(M_{2\ell}^+(u,i)\otimes L(v,j))&\cong \bigoplus_{\theta\in\mathfrak{I} }M_{2\ell}^+(u+v,i+j+\theta)\text{ and}\\\core(M_{2\ell}^-(u,i)\otimes L(v,j))&\cong \bigoplus_{\theta\in\mathfrak{I} }M_{2\ell}^-(u+v,i+j+\theta)
\end{align*} where $\I$ is defined in Proposition \ref{prop:simple tensor simple}.
\et

\bpf Using Propositions \ref{prop:case0} and \ref{prop:simple tensor simple} repeatedly, we have, on the one hand,
\begin{align*}
M_{2\ell}^+(0,0)\otimes L(u,i)\otimes L(v,j)&\cong M_{2\ell}^+(0,0)\otimes (L(u,i)\otimes L(v,j))\\
&\cong M_{2\ell}^+(0,0)\otimes \left(\bigoplus_{\theta\in\I}L(u+v,i+j+\theta)\oplus P_1\right)\\
&\cong \bigoplus_{\theta\in\I}M_{2\ell}^+(0,0)\otimes L(u+v,i+j+\theta) \oplus P_2\\
&\cong \bigoplus_{\theta\in\I}M_{2\ell}^+(u+v,i+j+\theta) \oplus P_3
\end{align*} and on the other hand
\begin{align*}
M_{2\ell}^+(0,0)\otimes L(u,i)\otimes L(v,j)&\cong (M_{2\ell}^+(0,0)\otimes L(u,i))\otimes L(v,j)\\
&\cong (M_{2\ell}^+(u,i)\oplus P_4)\otimes L(v,j)\\
&\cong (M_{2\ell}^+(u,i)\otimes L(v,j)) \oplus P_5
\end{align*} for some projective modules $P_1,\ldots,P_5$. The projective $P_5$ is necessarily a direct summand in $P_3$, therefore the result follows.

The proof of the second isomorphism is similar.
\epf

\section{Decomposition of the tensor product of two string modules of even length}\label{section:proof 4.22}

In this section, we determine the tensor products of two string modules of even length, in order to provide a correction to the statement and proof of \cite[Theorem 4.22]{EGST}.

We proved in \cite[Proposition 4.21]{EGST} that $M_2^+(u,i)\otimes M_2^-(v,j)$ is projective for any $(u,i)$ and $(v,j)$ in $\zz_n^2$. It then follows, using induction on $\ell$ and $t$ and the Auslander-Reiten sequences for the string modules of even length, that $M_{2\ell}^+(u,i)\otimes M_{2t}^-(v,j)$ is projective for any $(u,i)$ and $(v,j)$ in $\zz_n^2$ and any positive integers $\ell$ and $t.$

We denote by $\sHom_{\dland}(M,N)$ the space of stable homomorphisms from $M$ to $N$. We shall need the following isomorphisms.

\bl\label{lemma:stable adjoint} Let $U,$ $V$ and $W$ be any $\dland$-modules. Then there are isomorphisms
\begin{align*}
&\sHom_{\dland}(U\otimes V,W)\cong \sHom_{\dland}(U,V^*\otimes W)\\
&\sHom_{\dland}(U, V\otimes W)\cong \sHom_{\dland}(U\otimes V^*, W).
\end{align*}
\el

\bpf The adjoint functors $-\otimes V$ and $\Hom_k(V,-)$ take projective modules to projective modules, therefore by the Auslander-Kleiner theorem (see \cite[5.3.4]{Z}), there is an isomorphism $\Hom_{\dland}(U\otimes V,W)\cong \Hom_{\dland}(U,\Hom_k(V,W))$. Moreover, there is an isomorphism $\Hom_k(V,W)\cong V^*\otimes W$ of $\dland$-modules (see \cite[Section 4.2.C]{CP}), and the first isomorphism follows.

The second isomorphism is a consequence of the first isomorphism and of the isomorphism $V^{**}\cong V$ of $\dland$-modules (see \cite[Section 4.2.C]{CP}).
\epf

We now concentrate on tensor products of the form  $M_{2\ell}^+(u,i)\otimes M_{2t}^+(v,j)$ and   $M_{2\ell}^-(u,i)\otimes M_{2t}^-(v,j)$. In order to prove our results, we shall  need  the duals of some modules.

\bl\label{lemma:duals} Assume that $(u, i) \in \zz_n^2$ and
 $2i+u-1\not\equiv 0\pmod{d}$. Then  the following isomorphisms hold.
\begin{enumerate}[(i)]
\item $L(u,i)^*\cong L(-u,1-\sigma_u(i))$.
\item $M_{2\ell}^+(u,i)^*\cong M_{2\ell}^+(-u,1-i-\ell d).$
\item $M_{2\ell}^-(u,i)^*\cong M_{2\ell}^-(-u,1-i-(\ell-1)d).$
\end{enumerate}
\el

\bpf We start with the dual of a simple module. It is also a simple module, therefore there exists $(v,j)$ such that $L(u,i)^*\cong L(v,j).$ Set $N=\dim L(u,i)=\dim L(v,j)<d$. We must determine $v$ and $j$.

 The module $L(u,i)$ has a basis $\set{\varepsilon_t\,;\, 0\ppq t<N}$ given by $\varepsilon_t=X^t\tilde{F}_{u,i}$, using Proposition \ref{prop:simples} and the description of $L(u,i)$ preceding it. Let $\set{\varepsilon_t^*\,;\, 0\ppq t<N}$ be the dual basis. We have
\begin{itemize}
\item $(E_w\cdot \varepsilon_t^*)(\varepsilon_s)=\varepsilon_t^*(S(E_w)\varepsilon_s)=\varepsilon_t^*(E_{-w}\varepsilon_s)=
\begin{cases}
\varepsilon_t^*(\varepsilon_s)&\text{ if }-w=u\\0&\text{ otherwise}
\end{cases}$ \\so that $E_{-u}\cdot\varepsilon_t^*=\varepsilon_t^*$ and $E_{w}\cdot\varepsilon_t^*=0$ for $w\neq -u.$
By Proposition \ref{prop:simples} we must have $v=-u.$

\item Similarly, we have $X\cdot\varepsilon_t^*=
\begin{cases}
-q^{-i-t+1}\varepsilon_{t-1}^*&\text{ if }t>0\\0&\text{ if }t=0.
\end{cases}
$\\
In particular, $\varepsilon_0^*$ spans the kernel of the action of $X$. Then, by Proposition \ref{prop:simples}, the idempotent $e_j$ is the only one which does not
annihilate $\varepsilon_0^*$. We can easily check that $e_{1-\sigma_u(i)}\cdot\varepsilon_0^*= \varepsilon_0^*$, therefore $j=1-\sigma_u(i).$
\end{itemize}

\smallskip

We now consider the dual of an indecomposable string module $M$ of even length. Since $M\cong M^{**}$ as  $\dland$-modules, the dual module $M^*$ is also indecomposable.  Moreover, if $M$ has radical length $2$, then so has $M^*$ and $\Top (M^*)\cong (\soc M)^*$ and $\soc(M^*)\cong (\Top M)^*.$
We also note that if $x\in\zz_n$ with $x\not\equiv 0\pmod d$, then $\rep{x}+\rep{-x}=d.$

First consider $M_2^+(u,i)^*.$ We have $\soc(M_2^+(u,i)^*)\cong \Top (M_2^+(u,i))^*\cong L(u,i)^*\cong L(-u,1-\sigma_u(i))$ and $\Top(M_2^+(u,i)^*)\cong \soc (M_2^+(u,i))^*\cong L(u,\sigma_u(i))^*\cong L(-u,\sigma_u^2(i))=L(-u,1-i-d).$ Note that $\sigma_{-u}(1-i-d)=1-\sigma_u(i)$.
Moreover, $M_2^+(u,i)$ is not projective as a $\lndd$-module, therefore it has a direct summand of dimension strictly smaller than $d$. It follows that the $\lndd$-module  $M_2^+(u,i)^*$ also has  a direct summand of dimension strictly smaller than $d$, so that $M_2^+(u,i)^*$ is not projective as a $\lndd$-module. Therefore $M_2^+(u,i)^*\cong M_2^+(-u,1-i-d)$.

The dual $M_{2\ell}^+(u,i)^*$ is then obtained by induction, using the Auslander-Reiten sequences and the fact that dualising an  Auslander-Reiten sequence gives an Auslander-Reiten sequence.

\smallskip
The proof of (iii) is similar.
\epf

We shall now study the tensor product of string modules of even length, and we start with some special cases, which constitute the initial step in the proof of Theorem \ref{thm:even tensor even}.

\bp\label{prop:reduction} We have
\begin{align*}
\core(\evmp{t}{u}{i}\otimes \evmp{t}{v}{j})&\cong \bigoplus_{\theta\in\I} \left(\evmp{t}{u+v}{i+j+\theta}\oplus \evmp{t}{u+v}{\sigma_{u+v}^{2t-1}(u+v+\theta)}\right)\\
\core(\evmp{t}{u}{i}\otimes \evmp{(t+1)}{v}{j})&\cong \bigoplus_{\theta\in\I}\left( \evmp{t}{u+v}{i+j+\theta}\oplus \evmp{t}{u+v}{\sigma_{u+v}^{2t+1}(u+v+\theta)}\right)\\
\core(\evmpm{t}{u}{i}\otimes \evmpm{t}{v}{j})&\cong \bigoplus_{\theta\in\I} \left(\evmpm{t}{u+v}{i+j+\theta}\oplus \evmpm{t}{u+v}{\sigma_{u+v}^{-(2t-1)}(u+v+\theta)}\right)\\
\core(\evmpm{t}{u}{i}\otimes \evmpm{(t+1)}{v}{j})&\cong \bigoplus_{\theta\in\I}\left( \evmpm{t}{u+v}{i+j+\theta}\oplus \evmpm{t}{u+v}{\sigma_{u+v}^{-(2t+1)}(u+v+\theta)}\right)
\end{align*} where $\I$ is defined in Proposition \ref{prop:simple tensor simple}.
\ep

\bpf We start by proving the first two isomorphisms in the  case $(u,i)=(v,j)=(0,0)$, that is, \[\begin{cases}
\core(\evmp{t}{0}{0}\otimes \evmp{t}{0}{0})\cong \evmp{t}{0}{0}\oplus \evmp{t}{0}{\sigma_{0}^{2t-1}(0)}\\
\core(\evmp{t}{0}{0}\otimes \evmp{(t+1)}{0}{0})\cong  \evmp{t}{0}{0}\oplus \evmp{t}{0}{\sigma_{0}^{2t+1}(0)}
\end{cases}\]

Set $\ell=t$ or $\ell=t+1$.

We have an exact sequence
\[ 0\rightarrow L(0,1+(\ell-1)d)\rightarrow M_{2\ell}^+(0,0)\rightarrow \Omega^{\ell-1}(L(0,\sigma_0^{\ell-1}(0)))\rightarrow0. \] Tensoring on the left with $M_{2t}^+(0,0)$ gives an exact sequence
\[ 0\rightarrow M_{2t}^+(0,1+(\ell-1)d)\oplus P_1\rightarrow M_{2t}^+(0,0)\otimes M_{2\ell}^+(0,0)\rightarrow\Omega^{\ell-1}(M_{2t}^+(0,\sigma_0^{\ell-1}(0)))\oplus P_2\rightarrow 0 \] for some projective-injective modules $P_1$ and $P_2$. It follows that there is an exact sequence
\begin{equation}\label{eq:sequence tensor two even} 0\rightarrow M_{2t}^+(0,1+(\ell-1)d)\rightarrow \core(M_{2t}^+(0,0)\otimes M_{2\ell}^+(0,0))\oplus P_3\rightarrow M_{2t}^+(0,0)\rightarrow0 \end{equation} for some projective module $P_3$ (obtained by factoring out the split exact sequence $0\rightarrow P_1\rightarrow P_1\oplus P_2
\rightarrow P_2\rightarrow 0$). The exact sequence \eqref{eq:sequence tensor two even} is isomorphic to the following pullback
\[ \xymatrix{0\ar[r] &M_{2t}^+(0,1+(\ell-1)d)\ar[r]\ar@{=}[d] & E \ar[r]\ar[d]& M_{2t}^+(0,0)\ar[r]\ar[d]^\varphi&0\\
0\ar[r] &M_{2t}^+(0,1+(\ell-1)d)\ar[r] & P:=\bigoplus_{y=\ell}^{\ell+t-1}P(0,yd) \ar[r]^(.36)\pi&\Omega^{-1}( M_{2t}^+(0,1+(\ell-1)d))\cong M_{2t}^+(0,\ell d)\ar[r]&0} \] where $E=\set{(m,p)\in M_{2t}^+(0,0)\times P\,;\, \varphi(m)=\pi(p) }\cong \core(M_{2t}^+(0,0)\otimes M_{2\ell}^+(0,0))\oplus P.$

Assume for a contradiction that the sequence \eqref{eq:sequence tensor two even} is not split. Then $\varphi\neq0$. The modules $M_{2t}^+(0,0)$ and $M_{2t}^+(0,\ell d)$ have the same parity, therefore the image of $\varphi$ must contain a submodule of $M_{2t}^+(0,\ell d)$ of length $2,$ that is, $M_2^+(0,\ell d)$.
 The module $M_2^+(0,\ell d)$ has a generator of the form $\varphi(m)=\pi(e)$ for some $(m,e)\in M_{2t}^+(0,0)\times P(0,\ell d)\subset E.$ The submodule of $E$ generated by $(m,e)$ is projective, isomorphic to $P(0,\ell d)$. Since $P(0,\ell d)$ is injective, it is isomorphic to a direct summand of $E.$ Therefore $P(0,\ell d)$, which is a summand in $P$, is also a summand in $M_{2t}^+(0,0)\otimes M_{2\ell}^+(0,0)$.

Denote by $[M:L]$ the multiplicity of a simple module $L$ as a summand in a semisimple module $M.$ The simple module $L(0,\ell d)$ is in the socle of the projective summand $P(0,\ell d)$ of $M_{2t}^+(0,0)\otimes M_{2\ell}^+(0,0)$ (but does not occur in the socle of $M_{2t}^+(0,0)$), hence it follows from the above that
\[ [\soc (\core(M_{2t}^+(0,0)\otimes M_{2\ell}^+(0,0))) : L(0,\ell d)]<[\soc (M_{2t}^+(0,1+(\ell-1)d)):L(0,\ell d)].  \] We now compute these multiplicities. Set $t=r\frac{n}{d}+s$ with $0\ppq s<\frac{n}{d}.$ Recall that $\ell\pgq t.$

 We have
\begin{align*}
[\soc (M_{2t}^+(0,1+(\ell-1)d)):L(0,\ell d)]&=[\bigoplus_{x=\ell}^{\ell+t-1}L(0,xd):L(0,\ell d)]=[\bigoplus_{x=0}^{t-1}L(0,(x+\ell)d):L(0,\ell d)]\\&=\#\set{x\,;\, 0\ppq x\ppq t-1,\ x\equiv 0\pmod{\frac{n}{d}}}.
\end{align*}

 Using Lemmas \ref{lemma:stable adjoint} and \ref{lemma:duals} and Theorem \ref{thm:even tensor simple}, we have
\begin{align*}
[\soc (\core(M_{2t}^+(0,0)\otimes M_{2\ell}^+(0,0))) : L(0,\ell d)]&=\dim\sHom_{\dland}(L(0,\ell d),\core(M_{2t}^+(0,0)\otimes M_{2\ell}^+(0,0)))\\
&=\dim\sHom_{\dland}(L(0,\ell d),M_{2t}^+(0,0)\otimes M_{2\ell}^+(0,0))\\
&=\dim\sHom_{\dland}(L(0,\ell d)\otimes M_{2\ell}^+(0,0)^*,M_{2t}^+(0,0))\\
&=\dim\sHom_{\dland}(L(0,\ell d)\otimes M_{2\ell}^+(0,1-\ell d), M_{2t}^+(0,0))\\
&=\dim\sHom_{\dland}(M_{2\ell}^+(0,1), M_{2t}^+(0,0))\\
&=\dim\Ext^1_{\dland}(M_{2\ell}^+(0,d),M_{2t}^+(0,0))
\end{align*} which is known by Lemma \ref{lemma:dim Ext1}.

Finally,
\begin{align*}
\delta:=[\soc &(\core(M_{2t}^+(0,0)\otimes M_{2\ell}^+(0,0))) : L(0,\ell d)]-[\soc (M_{2t}^+(0,1+(\ell-1)d)) : L(0,\ell d)]\\
&=\#\set{y\,;\, 0\ppq y\ppq t-1,\ y\equiv 0\pmod{\frac{n}{d}}}-\#\set{y\,;\, 0\ppq y\ppq t-1,\ y\equiv \ell\pmod{\frac{n}{d}}}\\&\qquad +\#\set{y\,;\, 0\ppq y\ppq t-1,\ y\equiv t-\ell-1\pmod{\frac{n}{d}}}-\#\set{x\,;\, 0\ppq x\ppq t-1,\ x\equiv 0\pmod{\frac{n}{d}}}\\
&=\#\set{x\,;\, 0\ppq x\ppq t-1,\ x\equiv t-\ell-1\pmod{\frac{n}{d}}}-\#\set{x\,;\, 0\ppq x\ppq t-1,\ x\equiv \ell\pmod{\frac{n}{d}}}.
\end{align*}

We start with the case $\ell =t$. Then
\begin{align*}
0>\delta&=\#\set{x\,;\, 0\ppq x\ppq r\frac{n}{d}+s-1,\ x\equiv -1\pmod{\frac{n}{d}}}-\#\set{x\,;\, 0\ppq x\ppq r\frac{n}{d}+s-1,\ x\equiv s\pmod{\frac{n}{d}}}\\&=r-r=0
\end{align*} so that we have a contradiction.

In the case $\ell=t+1$,
\begin{align*}
0>\delta&=\#\set{x\,;\, 0\ppq x\ppq r\frac{n}{d}+s-1,\ x\equiv -2\pmod{\frac{n}{d}}}-\#\set{x\,;\, 0\ppq x\ppq r\frac{n}{d}+s-1,\ x\equiv s+1\pmod{\frac{n}{d}}}\\&=
\begin{cases}
(r+1)-(r+1)=0&\text{ if }s=\frac{n}{d}-1\\
r-r=0&\text{ if }s<\frac{n}{d}-1
\end{cases}
\end{align*} so that we have a contradiction.

Therefore the sequence \eqref{eq:sequence tensor two even} is  split in both cases and we have
\[\begin{cases}
\core(\evmp{t}{0}{0}\otimes \evmp{t}{0}{0})\cong \evmp{t}{0}{0}\oplus \evmp{t}{0}{\sigma_{0}^{2t-1}(0)}\quad\text{(when $\ell=t$)}\\
\core(\evmp{t}{0}{0}\otimes \evmp{(t+1)}{0}{0})\cong  \evmp{t}{0}{0}\oplus \evmp{t}{0}{\sigma_{0}^{2t+1}(0)}\quad\text{(when $\ell=t+1$)}.
\end{cases}\]

The result then follows, tensoring these isomorphisms with $L(u,i)$ on the left and $L(v,j)$ on the right, and using the commutativity of the tensor product and Proposition \ref{prop:simple tensor simple} on the tensor product of simple modules.

\medskip

The other two isomorphisms are proved in the same way: there is an exact sequence \[0\rightarrow M_{2t}^+(0,1-\ell d)\rightarrow \core(M_{2t}^-(0,0)\otimes M_{2\ell}^-(0,0))\oplus P_3\rightarrow M_{2t}^-(0,0)\rightarrow0 \] and if the middle term has a projective summand, then it must contain $P(0,-\ell d)$ as a summand. We then compare the multiplicities of the simple module $L(0,-\ell d)$ in $\soc(\core(M_{2t}^-(0,0)\otimes M_{2\ell}^-(0,0)))$ and in $\soc (M_{2t}^-(0,1-\ell d)).$
\epf

We may now rectify \cite[Theorem 4.22]{EGST}.

\bt\label{thm:even tensor even}  For any positive integer $t$ and any integer $\ell$ with $\ell\pgq t$, we have
\begin{align*}
\core(\evmp{\ell}{v}{j}\otimes\evmp{t}{u}{i})&\cong \bigoplus_{\theta\in\I} \left(\evmp{t}{u+v}{i+j+\theta}\oplus \evmp{t}{u+v}{\sigma^{2\ell-1}_{u+v}(i+j+\theta)}\right)\text{ and}\\
\core(\evmpm{\ell}{v}{j}\otimes\evmpm{t}{u}{i})&\cong \bigoplus_{\theta\in\I} \left(\evmpm{t}{u+v}{i+j+\theta}\oplus \evmpm{t}{u+v}{\sigma^{-(2\ell-1)}_{u+v}(i+j+\theta)}\right),
\end{align*}
 where $\I$ is defined in Proposition \ref{prop:simple tensor simple}.
\et

\bpf
 We prove the first isomorphism by  induction on $\ell,$ the proof of the second isomorphism is similar.

By Proposition \ref{prop:reduction}, we already know that the result is true for $\ell=t$ and for $\ell=t+1$.  Now take $\ell\pgq t+1$ and assume that the decomposition holds for any $\ell'$ with $t\ppq \ell'\ppq\ell$. We apply Lemma \ref{lemma:AC} with $M=\evmp{t}{u}{i}$ and $C=\evmp{\ell}{v}{j+d}.$  The module
$M^*\otimes M\otimes C\cong \evmp{t}{-u}{1-i-td}\otimes \evmp{t}{u}{i}\otimes \evmp{\ell}v{j+d}$ is the direct sum of a projective module and of indecomposable modules of length $2t<2\ell$ by the induction hypothesis and Proposition \ref{prop:reduction}, therefore $\evmp{\ell}{v}{j+d}$ is not a summand in $\evmp{t}{-u}{1-i-td}\otimes \evmp{t}{u}{i}\otimes \evmp{\ell}v{j+d}$.  It follows that $\A(\evmp{\ell}{v}{j+d})\otimes \evmp{t}{u}{i}$ splits, so that
\begin{align*}
&\core((\evmp{(\ell+1)}{v}{j}\otimes\evmp{t}{u}{i}) \oplus (\evmp{(\ell-1)}{v}{j+d}\otimes\evmp{t}{u}{i}))\\&\qquad\cong \core( (\evmp{\ell}{v}{j}\otimes\evmp{t}{u}{i}) \oplus (\evmp{\ell}{v}{j+d}\otimes\evmp{t}{u}{i})).
\end{align*}
Since $t\ppq \ell-1$, using the induction hypothesis we have
\begin{align*}
\core(\evmp{(\ell-1)}{v}{j+d}\otimes\evmp{t}{u}{i})&\cong \bigoplus_{\theta\in\I} \left(\evmp{t}{u+v}{i+j+\theta+d}\oplus \evmp{t}{u+v}{\sigma^{2\ell-1}_{u+v}(i+j+\theta)}\right)\\
\core(\evmp{\ell}{v}{j}\otimes\evmp{t}{u}{i})&\cong \bigoplus_{\theta\in\I} \left(\evmp{t}{u+v}{i+j+\theta}\oplus \evmp{t}{u+v}{\sigma^{2\ell-1}_{u+v}(i+j+\theta)}\right)\\
\core(\evmp{\ell}{v}{j+d}\otimes\evmp{t}{u}{i})&\cong\bigoplus_{\theta\in\I} \left( \evmp{t}{u+v}{i+j+\theta+d}\oplus \evmp{t}{u+v}{\sigma^{2\ell+1}_{u+v}(i+j+\theta)}\right)
\end{align*} therefore $\core(\evmp{(\ell+1)}{v}{j}\otimes\evmp{t}{u}{i})\cong \bigoplus_{\theta\in\I} \left(\evmp{t}{u+v}{i+j+\theta}\oplus \evmp{t}{u+v}{\sigma^{2\ell+1}_{u+v}(i+j+\theta)}\right)$ as required.

This concludes the induction.
\epf

\section{Decomposition of tensor products involving band modules}\label{section:band modules}

In \cite{EGST}, we were unable to give a general expression of the decomposition of the tensor product of a band module with another module. We can do this now, using the methods in the previous two sections.

We consider first the tensor product of a band module of shortest length with a simple module.

\bp (cf. \cite[Proposition 4.17]{EGST}) \label{prop:special band tensor L} For any $\lambda\in k\backslash\set{0}$, there exists $\mu\in k\backslash\set{0}$ such that $\core(C_\lambda^1(0,0)\otimes L(v,j))\cong C_\mu^1(v,j)$.
\ep

\bpf  The proof is similar to that of Proposition \ref{prop:special M2 tensor L}. Set $N=\dim L(v,j).$  We tensor the non-split exact sequence $0\rightarrow \Omega^{1-\frac{n}{d}}(L(0,\sigma_0^{\frac{n}{d}}(0)))\rightarrow C_\lambda^1(0,0)\rightarrow L(0,0)\rightarrow 0$ by $L(v,j)$. Using Proposition \ref{prop:simple tensor simple} as well as  $\dim L(0,0)=1$ and $\dim L(0,\sigma_0^{\frac{n}{d}}(0))\in\set{1,d-1}$ depending on the parity of $\frac{n}{d}$ (so that $\#\I=1$), we obtain  an exact sequence $0\rightarrow \Omega^{1-\frac{n}{d}}(L(v,\sigma_v^{\frac{n}{d}}(j)))\oplus P \rightarrow C_\lambda^1(0,0)\otimes L(v,j)\rightarrow L(v,j)\rightarrow 0$ with $P$ a projective module or $0$.
Since $P$ is injective, we have $C_\lambda^1(0,0)\otimes L(v,j)\cong U\oplus P$ for some module $U$ and we obtain an exact sequence $0\rightarrow\Omega^{1-\frac{n}{d}}(L(v,\sigma_v^{\frac{n}{d}}(j)))\rightarrow U\rightarrow   L(v,j)\rightarrow 0$.
The sequence cannot split by Lemma \ref{lemma:AC},
 so $U$ is an indecomposable module of length $2\frac{n}{d}$, such that $L(v,j)$ is contained in the top of $U$.
 Moreover, it follows from Property \ref{pty:charact band length pullback} and Remark \ref{rk:tensor restricted projective} that $C_\lambda^1(0,0)\otimes L(v,j)$, and hence $U$, is projective as a $\lnd$-module and as a $\lndd$-module, therefore  we must have $U=C_\mu^1(v,j)$ for some parameter $\mu\in k\backslash\set{0}.$
\epf

As in the case of string modules of even length, this proposition constitutes the initial step in the proof that $\core(C_\lambda^{\ell}(0,0)\otimes L(v,j))$ is a band module of the form $C_\mu^{\ell}(v,j).$ In order to continue, we shall need the extensions between band modules.

\bl\label{lemma:dim Ext1 band} Let $\ell$ and $t$ be positive integers. For any $(u,i)\in\zz_n^2$ with $2i+u-1\not\equiv 0\pmod d$ and any parameters $\lambda$ and $\mu$ in $k\setminus\set{0}$,
\[ \dim\Ext^1_{\dland}(C_\lambda^t(u,i),C_\mu^{\ell}(u,i))=
\begin{cases}
\min(t,\ell)&\text{ if }\lambda=\mu\\0&\text{ if }\lambda\neq \mu.
\end{cases}
 \]
Moreover, if $t=1$ and $\lambda=\mu,$ the equivalence class of the exact sequence
\[ 0\rightarrow C_\lambda^{\ell}(u,i)\rightarrow C_\lambda^{\ell+1}(u,i)\rightarrow C_\lambda^{1}(u,i)\rightarrow 0 \] is a basis of  $\Ext^1_{\dland}(C_\lambda^1(u,i),C_\mu^{\ell}(u,i))$
\el

\begin{proof} There is a parameter $\omega(\lambda)\in k\backslash\set{0}$ such that $\Omega(C_\lambda^t(u,i))=C_{\omega(\lambda)}^t(u,\sigma_u(i)).$
Apply $\Hom_{\dland}(-,C_\mu^{\ell}(u,i))$ to the exact sequence
\[ 0\rightarrow C_{\omega(\lambda)}^t(u,\sigma_u(i))\rightarrow \bigoplus_{r=1}^{\frac{n}{d}}P(u,i+rd)^t\rightarrow C_\lambda^t(u,i)\rightarrow 0 \] to get the exact sequence
\begin{align*}
0& \rightarrow \Hom_{\dland}(C_\lambda^t(u,i),C_\mu^{\ell}(u,i)) \rightarrow \bigoplus_{r=1}^{\frac{n}{d}}\Hom_{\dland}(P(u,i+rd),C_\mu^{\ell}(u,i))^t\\&\rightarrow \Hom_{\dland}(C_{\omega(\lambda)}^t(u,\sigma_u(i)),C_\mu^{\ell}(u,i)) \rightarrow \Ext^1_{\dland}(C_\lambda^t(u,i),C_\mu^{\ell}(u,i))\rightarrow 0.
\end{align*}

\begin{itemize}
\item Let $f$ be a non-zero map in $\Hom_{\dland}(C_{\omega(\lambda)}^t(u,\sigma_u(i)),C_\mu^{\ell}(u,i)) $. Since the modules have different parities, $f$ must map the top $\bigoplus_{x=0}^{\frac{n}{d}-1}L(u,\sigma_u^{2x+1}(i))^t$ of $C_{\omega(\lambda)}^t(u,\sigma_u(i))$ into the socle $\bigoplus_{x=0}^{\frac{n}{d}-1}L(u,\sigma_u^{2x+1}(i))^\ell$ of $C_\mu^{\ell}(u,i)$, therefore
$\dim \Hom_{\dland}(C_{\omega(\lambda)}^t(u,\sigma_u(i)),C_\mu^{\ell}(u,i))=t\ell\frac{n}{d}.$

\item $\dim\bigoplus_{r=1}^{\frac{n}{d}}\Hom_{\dland}(P(u,i+rd),C_\mu^{\ell}(u,i))^t=t\ell\frac{n}{d}$.
\item Therefore $\dim\Ext^1_{\dland}(C_\lambda^t(u,i),C_\mu^{\ell}(u,i))=\dim\Hom_{\dland}(C_\lambda^t(u,i),C_\mu^{\ell}(u,i))$ and we must prove that
\[ \dim\Hom_{\dland}(C_\lambda^t(u,i),C_\mu^{\ell}(u,i))=
\begin{cases}
0&\text{if }\lambda\neq \mu\\\min(t,\ell)&\text{if }\lambda= \mu.
\end{cases}
\]
Both $C_\lambda^t(u,i)$ and $C_\mu^{\ell}(u,i)$ are modules over the same block $\Bb_{u,i}$, therefore we may work over the basic algebra  $\Bb_{u,i}$ with band modules such that $b_0$ is the arrow acting as $J_t(\lambda)$ and $J_\ell(\mu)$ respectively. A morphism $f\in\Hom_{\dland}(C_\lambda^t(u,i),C_\mu^{\ell}(u,i))$ is then completely determined by a linear map $f_0\in\Hom_k(k^t,k^\ell)$ such that $J_\ell(\mu)f_0=f_0J_t(\lambda)$. Such a map is zero if $\lambda\neq \mu$, and if $\lambda=\mu$ its matrix is of the form $\begin{pmatrix}A\\0_{\ell-t}\end{pmatrix}$ if $t\ppq \ell$ and $\begin{pmatrix}A & 0_{t-\ell}\end{pmatrix}$ if $t\pgq \ell$, where $A$ is an upper triangular Toeplitz matrix, that is, of the form $A=\begin{pmatrix}a_1&a_2&\cdots&a_m\\0&a_1&\ddots&\vdots\\\vdots&&\ddots&a_2\\0&\cdots&\cdots&a_1 \end{pmatrix}$, with $m=\min(\ell,t)$ and $a_1,\ldots,a_m$ in $k$.
\qedhere
\end{itemize}
\end{proof}

We can now determine the tensor product of a longer band module with a simple module.

\bp\label{prop:case0 band} Let $\ell$ be a positive integer and let $\lambda$ be in $k\backslash\set{0}$. Let $\mu\in k\backslash\set{0}$ be as in Proposition \ref{prop:special band tensor L}. Then
\[ \core(C_\lambda^{\ell}(0,0)\otimes L(v,j))\cong C_\mu^{\ell}(v,j). \]
\ep

\begin{proof}
The proof is by induction on $\ell.$ Proposition \ref{prop:special band tensor L} shows that the result is true for $\ell=1.$

 Now assume that $\core(C_\lambda^{\ell}(0,0)\otimes L(v,j))\cong C_\mu^{\ell}(v,j)$ for a given $\ell\pgq 1.$ There is an exact sequence
\[ 0\rightarrow C_\lambda^{\ell}(0,0)\rightarrow C_\lambda^{\ell+1}(0,0)\rightarrow C_\lambda^1(0,0)\rightarrow0. \]  Tensoring with $L(v,j)$ gives an exact sequence
\[ 0\rightarrow C_\mu^{\ell}(v,j)\oplus P_1\rightarrow C_\lambda^{\ell+1}(0,0)\otimes L(v,j)\rightarrow C_\mu^1(v,j) \oplus P_2\rightarrow 0  \] with $P_1$ and $P_2$ projective-injective modules, so that, factoring out the split exact sequence $0\rightarrow P_1\rightarrow P_1\oplus P_2\rightarrow P_2\rightarrow 0$, we have an exact sequence
\[ 0\rightarrow C_\mu^{\ell}(v,j)\rightarrow \core(C_\lambda^{\ell+1}(0,0)\otimes L(v,j))\oplus P\rightarrow C_\mu^1(v,j) \rightarrow 0  \] for some projective module $P.$ Moreover, since $L(v,j)$ is a splitting trace module, by Lemma \ref{lemma:AC} this sequence is not split. By Lemma \ref{lemma:dim Ext1 band}, we have $P=0$ and $\core(C_\lambda^{\ell+1}(0,0)\otimes L(v,j))=C_\mu^{\ell+1}(v,j)$ thus proving the induction.
\end{proof}

As we mentioned in Section \ref{section:list modules},  the parameter $\lambda$ of the module $C_\lambda^\ell(u,i)$ is not  well defined, due to the fact that defining it requires a fixed labelling of the vertices and  arrows of the basic algebra of $\dland$, which we do not have. Moreover, when computing tensor products, we work with $\dland$-modules (and not over the basic algebra).  The method in \cite[Example 4.23]{EGST} allows us to determine the parameter $\mu$ in specific examples, once the labelling choices are made, but it does not give a general rule.

However, Proposition \ref{prop:case0 band} allows us to fix the parameters coherently once and for all in the following way.

\medskip

\paragraph{\textbf{Convention on parameters.}} Define $C_\lambda^\ell(0,0)$ as in Subsection \ref{subsec:even}, with the vertex $\epsilon_0$ of $\Bb_{0,0}$ corresponding to the trivial module $L(0,0)$. We then set \[C_\lambda^{\ell}(v,j)=\core(C_\lambda^{\ell}(0,0)\otimes L(v,j)).\]

With the definitions recalled in Subsection \ref{subsec:even} we have $C_\lambda^\ell(u,i)\cong C_\lambda^\ell(u,\sigma_u^2(i))=C_\lambda^\ell(u,i+d)=\Omega^2(C_\lambda^\ell(u,i))$, and this is compatible with the convention above.

\medskip

The tensor product of any band module with a simple module is an immediate consequence of this convention.

\bp\label{prop:band tensor simple}  Fix an integer $\ell\pgq 1 $ and a parameter $\lambda\in k\backslash\set{0}$. For any $(u,i)$ and $(v,j)$ in $\zz_n^2$ with $2i+u-1\not\equiv 0\pmod{d}$ and $2j+v-1\not\equiv 0\pmod d$, we have
\[ \core(C_\lambda^\ell(u,i)\otimes L(v,j))\cong \bigoplus_{\theta\in\mathfrak{I} }C_\lambda^{\ell}(u+v,i+j+\theta) \]
 where $\I$ is defined in Proposition \ref{prop:simple tensor simple}.
\ep

\br As a consequence, we know the tensor product of band modules with string modules of odd length up to projectives.
\er

We can now consider the tensor product of band modules with modules of even length. The first result concerns tensor products of band modules with different parameters and  tensor products of band modules with string modules of even length.

\bp\label{prop:tensor different bands} Let $(u,i)$ and $(v,j)$ be in $\zz_n^2$ such that $d$ does not divide $2i+u-1$ or $2j+v-1$, let  $\ell$ and $t$ be any positive integers.

For any $\lambda,\mu$ in $k\backslash\set{0}$ such that $\lambda\neq \mu$, the module $C_\lambda^{\ell}(u,i)\ot C_\mu^{t}(v,j)$ is projective.

For any $\lambda$ in $k\backslash\set{0}$, the module $C_\lambda^{\ell}(u,i)\ot M_{2t}^{\pm}(v,j)$ is projective.

\ep

\bpf We start by proving that $C_\lambda^{1}(0,0)\ot C_\mu^{1}(0,0)$ is projective. First note that  for any parameter $\lambda\in k\backslash\set{0}$  there exist parameters $\omega(\lambda),\omega'(\lambda)\in k\backslash\set{0}$ such that $\Omega^{-1}(C_\lambda^1(0,0))=C_{\omega(\lambda)}^1(0,1)$ and $\Omega^{-1}(C_\lambda^1(0,1))=C_{\omega'(\lambda)}^1(0,0)$, and moreover that $\omega'(\omega(\lambda))=\lambda$ so that $\omega$ is an injective map.

Assume first that $\frac{n}{d}$ is odd.
\begin{itemize}
\item There is an exact sequence
$0\rightarrow L(0, \sigma_0^{\frac{n}{d}}(0))\rightarrow \Omega^{-\frac{n}{d}}(L(0,0))
\rightarrow C_\lambda^1(0,0)\rightarrow 0$. We tensor this
with $C_\mu^1(0,0)$ on the right and factor out projective-injectives.
Using the fact that $C_{\mu}^1(0, \sigma_0^{\frac{n}{d}}(0))\cong C_{\mu}^1(0,1)$,  the resulting exact sequence is

\[ 0\rightarrow C_\mu^1(0, 1) \xrightarrow[(f_1,f_2)]{} \Omega^{-\frac{n}{d}}(C_\mu^1(0,0))\oplus \text{ projective}\cong C_{\omega(\mu)}^1(0,1)\oplus P\rightarrow C_\lambda^1(0,0)\ot C_\mu^1(0,0)\rightarrow 0 .\] If $f_1\neq0$, then $\omega(\mu)=\mu$ and $f_1$ is an isomorphism, therefore  $C_\lambda^1(0,0)\ot C_\mu^1(0,0)\cong P$. If $f_1=0,$ then $C_\lambda^1(0,0)\ot C_\mu^1(0,0)\cong C_{\omega(\mu)}^1(0,1)\oplus \Omega^{-1}(C_{\mu}^1(0,1)\oplus \text{ projective} $.

So either  $C_\lambda^1(0,0)\ot C_\mu^1(0,0)$ is projective or $C_\lambda^1(0,0)\ot C_\mu^1(0,0)\cong C_{\omega(\mu)}^1(0,1)\oplus C_{\omega'(\mu)}^1(0,0)\oplus \text{ projective} $.
\item Exchanging $\lambda$ and $\mu$ shows that  either  $C_\lambda^1(0,0)\ot C_\mu^1(0,0)$ is projective or $C_\lambda^1(0,0)\ot C_\mu^1(0,0)\cong C_{\omega(\lambda)}^1(0,1)\oplus C_{\omega'(\lambda)}^1(0,0)\oplus \text{ projective} $.

\item Therefore, since $\omega$ is injective,  either $\lambda=\mu$ or $C_\lambda^1(0,0)\ot C_\mu^1(0,0)$ is projective.
\end{itemize}

If $\frac{n}{d}$ is even,  the same arguments work if we start with the
exact sequence
\[0\to C_{\lambda}^1(0,0)\to \Omega^{\frac{n}{d}}(L(0,0))\to L(0, \sigma_0^{\frac{n}{d}}(0))\to 0\]
The proof that $C_\lambda^{1}(0,0)\ot M_{2}^{\pm}(0,xd)$ for $x\in\zz$ is projective is similar, using the exact sequences $0\rightarrow  L(0,1+(x-1)d)\rightarrow  \Omega^{-1}(L(0,xd))\rightarrow M_2^+(0,xd)\rightarrow  0$ and $0\rightarrow  L(0,1+xd)\rightarrow  \Omega^{-1}(L(0,xd))\rightarrow M_2^-(0,xd)\rightarrow  0$. Then, using induction, the Auslander-Reiten sequences for string and band modules, and the fact that the tensor product of any module by a projective module is projective, it follows that $C_\lambda^{\ell}(0,0)\ot C_\mu^{t}(0,0)$ for $\lambda\neq \mu$ and $C_\lambda^{\ell}(0,0)\ot M_{2t}^\pm(0,xd)$  are projective.

The general case follows by tensoring with $L(u,i)\otimes L(v,j)$.
\epf

As a consequence, we can give  the dual of a band module.

\bc\label{cor:dual band} The dual of $C_\lambda^\ell(u,i)$ is  $C_\lambda^\ell(-u,1-i)$.
\ec

\bpf First, as in the case of string modules of even length, and using the fact that the dual of a projective $\lnd$-module (respectively $\lndd$-module) is projective because both algebras are self-injective as finite dimensional Hopf algebras, it can be shown that there exists $\mu\in k\backslash\set{0}$ such that $C_\lambda^\ell(u,i)^*\cong C_\mu^\ell(-u,1-i)$. Moreover, we have
\begin{align*}
0\neq &\sHom_{\dland}(C_\lambda^\ell(u,i),C_\lambda^\ell(u,i))\cong \sHom_{\dland}(C_\lambda^\ell(u,i)\otimes C_\lambda^\ell(u,i)^*,L(0,0))\\&\cong \sHom_{\dland}(C_\lambda^\ell(u,i)\otimes C_\mu^\ell(-u,1-i),L(0,0))
\end{align*}
 therefore $C_\lambda^\ell(u,i)\otimes C_\mu^\ell(-u,1-i)$ is not projective and by Proposition \ref{prop:tensor different bands} it follows that $\mu=\lambda$.
\epf

Another consequence is that we can determine the syzygy of any band module.

\bc \label{cor:syzygies} Let  $(u,i)$ be in $\zz^2_n$ with $2i+u-1\not\equiv 0\pmod d$ and let $\lambda\in k\backslash\set{0}$ be a parameter.  Then for any positive integer $\ell$ we have
\[ \Omega(C_\lambda^{\ell}(u,i))\cong C_{\lambda}^{\ell}(u,\sigma_u(i)) \cong \Omega^{-1}(C_\lambda^{\ell}(u,i)). \]
\ec

\bpf We start with the syzygy of $\band\lambda\ell00$.  We know that  there exists a parameter $\mu\in k\backslash\set{0}$ such that $\Omega(C_\lambda^\ell(0,0))=C_{\mu}^\ell(0,1)$. By Lemma \ref{lemma:dim Ext1 band}, the space $\Ext^1_{\dland}(C_{\mu}^\ell(0,1),C_{\mu}^\ell(0,1))$ is not zero. Moreover, there are isomorphisms
\begin{align*}
\Ext^1_{\dland}(C_{\mu}^\ell(0,1),C_{\mu}^\ell(0,1))&\cong \sHom_{\dland}(\Omega(\band{\mu}\ell01),\band{\mu}\ell01\otimes L(0,0))\\
&\cong \sHom_{\dland}(\Omega^2(\band{\lambda}\ell00)\otimes \band{\mu}\ell01^*,L(0,0))\\
&\cong\sHom_{\dland}(\band{\lambda}\ell00\otimes\band{\mu}\ell00,L(0,0))
\end{align*} and therefore $\band{\lambda}\ell00\otimes\band{\mu}\ell00$ is not projective. It follows from  Proposition \ref{prop:tensor different bands} that $\mu=\lambda$ and therefore that $\Omega(C_\lambda^\ell(0,0))=C_{\lambda}^\ell(0,1)$.

We then have the following isomorphisms.
\begin{align*}
\Omega(C_\lambda^\ell(u,i))&\cong\Omega(\core(C_\lambda^\ell(0,0)\otimes L(u,i)))\cong\core(\Omega(C_\lambda^\ell(0,0)\otimes L(u,i)))\\
&\cong \core(\Omega(C_\lambda^\ell(0,0))\otimes L(u,i))\cong\core(C_{\lambda}^\ell(0,1)\otimes L(u,i))
\\
&\cong\core(C_{\lambda}^\ell(0,0)\otimes L(0,1)\otimes L(u,i))\cong \core(C_{\lambda}^\ell(0,0)\otimes L(u,\sigma_u(i)))\cong C_{\lambda}^{\ell}(u,\sigma_u(i)).
\end{align*}

Finally, since $\Omega^{2}(C_\lambda^\ell(u,i))\cong C_\lambda^\ell(u,i)$, we also have $\Omega^{-1}(C_\lambda^{\ell}(u,i))\cong\Omega(C_\lambda^\ell(u,i))$.
\epf

The next result gives the tensor product of any two band modules with the same parameter up to projectives.

\bt\label{thm:band tensor band}  For any positive integer $t$ and any integer $\ell$ with $\ell\pgq t$, we have \[
\core(\band{\lambda}{\ell}{v}{j}\otimes\band{\lambda}{t}{u}{i})\cong \bigoplus_{\theta\in\I} \left(\band{\lambda}{t}{u+v}{i+j+\theta}\oplus \band{\lambda}{t}{u+v}{\sigma_{u+v}(i+j+\theta)}\right).
\] where  $\I$ is defined  in Proposition \ref{prop:simple tensor simple}.
\et

\bpf We prove the  isomorphism in the  case $(u,i)=(v,j)=(0,0)$, that is,
\[ \core(\band\lambda{\ell}{0}{0}\otimes\band{\lambda}{t}{0}{0})\cong \band{\lambda}{t}{0}{0}\oplus \band{\lambda}{t}{0}{1}
 \] and the result follows by taking the tensor product with $L(u,i)\otimes L(v,j).$

By working with modules over the basic algebra,
we can see that there is an exact sequence
\[ 0\rightarrow L(0,1-d)\rightarrow \band{\lambda}{\ell}{0}{0}\rightarrow \Omega^{\ell\frac{n}{d}-1}(L(0,\sigma_0^{\ell\frac{n}{d}-1}(0)))\rightarrow0. \] Tensoring on the right with $\band{\lambda}{t}{0}{0}$ gives an exact sequence
\[ 0\rightarrow \band{\lambda}{t}{0}{1}\oplus P_1\rightarrow \band{\lambda}{\ell}{0}{0}\otimes \band{\lambda}{t}{0}{0}\rightarrow\Omega^{\ell\frac{n}{d}-1}(\band{\lambda}{t}{0}{\sigma_0^{\ell\frac{n}{d}-1}(0)})\oplus P_2\rightarrow 0 \] for some projective-injective modules $P_1$ and $P_2$. It follows that there is an exact sequence
\begin{equation}\label{eq:sequence tensor two band} 0\rightarrow \band{\lambda}{t}{0}{1}\rightarrow \core(\band{\lambda}{\ell}{0}{0}\otimes \band{\lambda}{t}{0}{0})\oplus P_3\rightarrow \band{\lambda}{t}{0}{0}\rightarrow0 \end{equation} for some projective module $P_3$ (obtained by factoring out the split exact sequence $0\rightarrow P_1\rightarrow P_1\oplus P_2
\rightarrow P_2\rightarrow 0$). The exact sequence \eqref{eq:sequence tensor two band} is isomorphic to the following pullback
\[ \xymatrix{0\ar[r] &\band{\lambda}{t}{0}{1}\ar[r]\ar@{=}[d] & E \ar[r]\ar[d]& \band{\lambda}{t}{0}{0}\ar[r]\ar[d]^\varphi&0\\
0\ar[r] &\band{\lambda}{t}{0}{1}\ar[r] & P:=\bigoplus_{y=0}^{\frac{n}{d}-1}P(0,yd)^t \ar[r]^(.45)\pi&\Omega^{-1}( \band{\lambda}{t}{0}{1})\cong \band{\lambda}{t}00\ar[r]&0} \] where $E=\set{(m,p)\in \band{\lambda}{t}{0}{0}\times P\,;\, \varphi(m)=\pi(p) }\cong \core(\band{\lambda}{\ell}{0}{0}\otimes \band{\lambda}{t}{0}{0})\oplus P.$

Assume for a contradiction that the sequence \eqref{eq:sequence tensor two band} is not split. Then $\varphi\neq 0$ so that $\varphi$ is a non-zero endomorphism of $C_\lambda^t(0,0)$.  Let $m\in C_{\lambda}^t(0,0)$ be an element that is
not in the radical and such that $\varphi(m)\neq 0$, we may assume that
the submodule generated by $m$ has a simple top. Then $\varphi(m) = \pi(e)$ where
$e\in P$ and $e$ generates an indecomposable projective summand.
 The element
$(m, e)$ belongs to $E$ and it generates an indecomposable
projective module. It is
then a summand of $E$ since projectives are injective.
Let this summand be $P(0, yd)$ for some $y$, it
 is then also a summand in
 $\band{\lambda}{\ell}{0}{0}\otimes \band{\lambda}{t}{0}{0}$.
Denote by $[M:L]$ the multiplicity of a simple module $L$ as a summand in a semisimple module $M.$ The simple module $L(0, yd)$ is in the socle of the projective summand $P(0,yd)$ of $\band{\lambda}{\ell}{0}{0}\otimes \band{\lambda}{t}{0}{0}$ (but does not occur in the socle of $C_\lambda^t(0,0)$), hence it follows from the above that
\[ [\soc (\core(\band{\lambda}{\ell}{0}{0}\otimes \band{\lambda}{t}{0}{0})) : L(0, yd)]<[\soc C_\lambda^t(0,1):L(0,yd)]=t.  \]
 Using Lemma \ref{lemma:stable adjoint},
Proposition \ref{prop:case0 band}, Corollary \ref{cor:dual band} and the fact that $C_{\lambda}^t(0, 1+yd)\cong C_{\lambda}^t(0, 1) $, we have
\begin{align*}
[\soc (\core(\band{\lambda}{\ell}{0}{0}\otimes \band{\lambda}{t}{0}{0})) : L(0,yd)]&=\dim\sHom_{\dland}(L(0,yd),\core(\band{\lambda}{\ell}{0}{0}\otimes \band{\lambda}{t}{0}{0}))\\
&=\dim\sHom_{\dland}(L(0,yd),\band{\lambda}{\ell}{0}{0}\otimes \band{\lambda}{t}{0}{0})\\
&=\dim\sHom_{\dland}(L(0, yd)\otimes C_\lambda^t(0,0)^*,\band{\lambda}{\ell}{0}{0})\\
&=\dim\sHom_{\dland}(\band{\lambda}{t}{0}{1},\band{\lambda}{\ell}{0}{0})\\
&=\dim\Ext^1_{\dland}(C_{\lambda}^t(0,0),\band{\lambda}{\ell}{0}{0})=t
\end{align*}
by Lemma \ref{lemma:dim Ext1 band} since $t\ppq \ell$.
 Therefore we have obtained a contradiction, the sequence \eqref{eq:sequence tensor two band} splits, and the result follows.
\epf

 In many of our proofs, we have worked with $\Bb_{0,0}$-modules, then tensored with non-projective simple modules to obtain the general result we were seeking. This approach can be formalised as follows.

\bt\label{thm:stable equiv blocks} Let $L(u,i)$ be a non-projective simple module. Then $-\otimes L(u,i)$ induces a stable equivalence between $\Bb_{0,0}$ and $\Bb_{u,i}$.
\et

The proof uses the following lemma.

\bl\label{lemma:tensor dual distinct blocks} Let $L(u,i)$ be a non-projective simple module of dimension $N.$ There is an isomorphism \[\core(L(u,i)\otimes L(u,i)^*)\cong \bigoplus _{\tau=1-\min(N,d-N)}^0L(0,\tau),\] and the blocks $\Bb_{0,\tau}$, with $1-\min(N,d-N)\ppq \tau\ppq 0$, are pairwise distinct.
\el

\bpf The isomorphism follows from Lemma \ref{lemma:duals} which states that $L(u,i)^*\cong L(-u,1-\sigma_u(i))$, Proposition \ref{prop:simple tensor simple}, and the fact that $1-\sigma_u(i)+i=1-N$. Since $1-\min(N,d-N)>-\frac{d}{2}$, we need only prove that the blocks $\Bb_{0,\tau}$ with $-\frac{d}{2}<\tau\ppq 0$ are pairwise different.

The block $\Bb_{0,\tau}$ contains precisely the simple modules $L(0,j)$ with $j$ in the $\sigma_0$-orbit of $\tau$, that is, $j\in\set{\tau+td,\sigma_0(\tau)+td\,;\, 0\ppq t<\frac{n}{d}}$ (recall that $j$ is taken modulo $n$). There are precisely two representatives of the $\sigma_0$-orbit of $\tau$ in $]-d,0]$.
 Moreover, if $-d<j\ppq -\frac{d}{2}$, then $-1-2d<2j-1\ppq -d-1$ so $\sigma_0(j)=d+j-\rep{2j-1}=d+j-(2j-1+2d)=-j-d+1$, and $-\frac{d}{2}+1\ppq -j-d+1<1$ therefore $-\frac{d}{2}<\sigma_0(j)\ppq 0$. Similarly, if $-\frac{d}{2}<j\ppq 0$ then $-d<\sigma_0^{-1}(j)=\sigma_0(j)-d\ppq -\frac{d}{2}$. It follows that the $\sigma_0$-orbit of $\tau$ has precisely one representative in $\mathopen{]}-\frac{d}{2}\mathclose{}\mathpunct{},0\mathclose{]}$, which proves our claim.
\epf

\bpf[Proof of Theorem \ref{thm:stable equiv blocks}]
Define a functor $F:\Bb_{0,0}\text{-}\smod\rightarrow \dland\text{-}\smod$ by $F(M)=\core(M\otimes L(u,i))$ on objects and $F(f)=f\otimes \id_{L(u,i)}$ restricted and co-restricted to the cores of the modules on morphisms.

We first determine the image of a simple module under $F$. The simple modules in $\Bb_{0,0}$ are the modules $L(0,td)$, which have dimension $1,$ and the modules $L(0,1+td)$, which have dimension $d-1$, for $0\ppq t<\frac{n}{d}$. By Proposition \ref{prop:simple tensor simple}, we have
\begin{align*}
F(L(0,td))&=\core(L(0,td)\otimes L(u,i))\cong L(u,i+td)=L(u,\sigma_u^{2t}(i))\\
F(L(0,1+td))&=\core(L(0,1+td)\otimes L(u,i))\cong L(u,\sigma_u(i)+td)=L(u,\sigma_u^{2t+1}(i))
\end{align*} which are simple $\Bb_{u,i}$-modules.
By induction on the length of modules, it follows that $F$ sends any non-projective $\Bb_{0,0}$-module to a non-projective $\Bb_{u,i}$-module. Therefore $F$ induces a functor $\Bb_{0,0}\text{-}\smod\rightarrow \Bb_{u,i}\text{-}\smod$ which we denote also by $F$.

We can be more specific. We have the following:
\begin{align*}
F(M_{2\ell}^\pm(0,0))&=\core(M_{2\ell}^\pm(0,0)\otimes L(u,i))\cong M_{2\ell}^\pm(u,i)\\
F(C_\lambda^\ell(0,0))&=\core(C_\lambda^\ell(0,0)\otimes L(u,i))\cong C_\lambda^\ell(u,i).
\end{align*}
Now note that  the functor $-\otimes L(u,i)$ defined on $\dland$-mod takes projectives to projectives and is exact, therefore it commutes with $\Omega$. It follows that $F$ commutes with $\Omega.$ Therefore we can give the image by $F$ of any non-projective indecomposable module:
\begin{equation*}\label{images by F}\begin{aligned}
F(\Omega^m(L(0,\sigma_0^t(0))))&\cong \Omega^m(F(L(0,\sigma_0^t(0))))\cong \Omega^m(L(u,\sigma_u^t(i))),\\
F(M_{2\ell}^\pm(0,\sigma_0^t(0)))&\cong F(\Omega^{-(\pm t)}(M_{2\ell}^\pm(0,0)))\cong \Omega^{-(\pm
  t)}(F(M_{2\ell}^\pm(0,0)))\cong \Omega^{-(\pm t)}(M_{2\ell}^\pm(u,i))\cong M_{2\ell}^\pm(u,\sigma_u^t(i)),\\
F(C_\lambda^\ell(0,\sigma_0^t(0)))&\cong F(\Omega^{-t}(C_\lambda^\ell(0,0)))\cong \Omega^{-t}(F(C_\lambda^\ell(0,0)))\cong \Omega^{-t}(C_\lambda^\ell(u,i))\cong C_\lambda^\ell(u,\sigma_u^t(i)).
\end{aligned}\end{equation*}

Consequently, for any $M'\in\Bb_{u,i}\text{-}\smod$, there exists $M\in \Bb_{0,0}\text{-}\smod$ such that $M'\cong F(M)$.

Finally, in order to prove that $F$ is indeed an equivalence of categories, we must prove that for any modules $M_1$ and $M_2$ in $\Bb_{0,0}$-\smod, we have $\sHom_{\Bb_{u,i}}(F(M_1),F(M_2))\cong \sHom_{\Bb_{0,0}}(M_1,M_2).$ We have
\begin{align*}
\sHom_{\Bb_{u,i}}(F(M_1),F(M_2))&\cong \sHom_{\dland}(M_1\otimes L(u,i),M_2\otimes L(u,i))\\
&\cong \sHom_{\dland}(M_1\otimes L(u,i)\otimes L(u,i)^*,M_2)\text{ by Lemma \ref{lemma:stable adjoint}}\\
&\cong \bigoplus_{\tau=1-\min(N,1-N)}^0\sHom_{\dland}(M_1\otimes L(0,\tau),M_2)\text{ by Lemma \ref{lemma:tensor dual distinct blocks}}\\
&\cong \bigoplus_{\tau=1-\min(N,1-N)}^0\sHom_{\Bb_{0,0}}(M_1\otimes L(0,\tau),M_2)
\end{align*} since $M_2$ belongs to $\Bb_{0,0}$. However, by Lemma \ref{lemma:tensor dual distinct blocks} and the beginning of this proof, the modules $M_1\otimes L(0,\tau)$ belong to pairwise distinct blocks, and the only one belonging to $\Bb_{0,0}$ is $M_1\otimes L(0,0)\cong M_1$. Therefore $\sHom_{\Bb_{u,i}}(F(M_1),F(M_2))\cong \sHom_{\Bb_{0,0}}(M_1,M_2)$.
\epf

\section{Description of the stable Green ring of $\dland$}\label{sec:summary}

Combining the results in this paper with the results in \cite{EGST}, we now have a complete description of the stable Green ring of $\dland$, which we give in Table \ref{table:summary} on page \pageref{table:summary}.

The stable Green ring is commutative so that we may assume for instance that  $\ell\pgq t$ in Table \ref{table:summary}.

\afterpage{\clearpage
  \begin{landscape}
\begin{center}
  \vfill
\def\arraystretch{2}
  \begin{tabular}{|c||c|c|c|c|}
    \hline
 $\otimes$   &$\Omega^m(L(v,j))$&$M_{2\ell}^+(v,j)$&$M_{2\ell}^-(v,j)$&$C_\mu^\ell(v,j)$\\\hline\hline
    $\Omega^n(L(u,i))$&$\displaystyle\bigoplus_{\theta\in\I}\Omega^{m+n}(L(w,i+j+\theta))$&$\displaystyle\bigoplus_{\theta\in\I}M^+_{2\ell}(w,\sigma_{w}^{-n}(i+j+\theta))$&$\displaystyle\bigoplus_{\theta\in\I}M^-_{2\ell}(w,\sigma_{w}^{n}(i+j+\theta))$&
\begin{tabular}{l}
if $n$ is even: $   \displaystyle\bigoplus_{\theta\in\I}C_\mu^\ell(w,i+j+\theta)$\\
if $n$ is odd:     $ \displaystyle\bigoplus_{\theta\in\I}C_\mu^\ell(w,\sigma_{w}(i+j+\theta))$
\end{tabular}
\\
&\cite[Theorem 4.1 and \S\ 4.2]{EGST}& Theorem \ref{thm:even tensor simple} & Theorem \ref{thm:even tensor simple}
                                                                      & Proposition \ref{prop:band tensor simple}
\\\hline
    $M_{2t}^+(u,i)$& $\displaystyle\bigoplus_{\theta\in\I}M^+_{2t}(w,\sigma_{w}^{-m}(i+j+\theta))$&
\begin{tabular}{c}
$\displaystyle\bigoplus_{\theta\in\I} \Big(\evmp{t}{w}{i+j+\theta}\oplus$\\\quad$\evmp{t}{w}{\sigma^{2\ell-1}_{w}(i+j+\theta)}\Big) $
\end{tabular}
&0&0\\& Theorem \ref{thm:even tensor simple} &  Theorem \ref{thm:even tensor even} & \cite[Proposition
                                                                               4.21]{EGST} &
                                                                                             Proposition \ref{prop:tensor different bands}
\\\hline
    $M_{2t}^-(u,i)$&$\displaystyle\bigoplus_{\theta\in\I}M^-_{2t}(w,\sigma_{w}^{m}(i+j+\theta))$ &0&
\begin{tabular}{c}
$\displaystyle\bigoplus_{\theta\in\I} \Big(\evmpm{t}{w}{i+j+\theta}\oplus$\\\quad$\evmpm{t}{w}{\sigma^{-(2\ell-1)}_{w}(i+j+\theta)}\Big)$
\end{tabular}
&$0$\\
& Theorem \ref{thm:even tensor simple} &  \cite[Proposition
                                                                               4.21]{EGST} &Theorem \ref{thm:even tensor even} & 
                                                                                             Proposition \ref{prop:tensor different bands}
\\
\hline
    $C_\lambda^t(u,i)$&
\begin{tabular}{l}
if $m$ is even: $\displaystyle\bigoplus_{\theta\in\I}C_\lambda^t(w,i+j+\theta)$\\
if $m$ is odd:    $\displaystyle\bigoplus_{\theta\in\I}C_\lambda^t(w,\sigma_{w}(i+j+\theta))$
\end{tabular} 
 &0&0&
\begin{tabular}{l}
 if $\lambda=\mu$: \quad $\displaystyle\bigoplus_{\theta\in\I} \Big(\band{\lambda}{t}{w}{i+j+\theta}\oplus$\\
      \hphantom{if $\lambda=\mu$:
  \quad}$\band{\lambda}{t}{w}{\sigma_{w}(i+j+\theta)}\Big)$\\otherwise: \quad $0$
\end{tabular}
 \\& Proposition \ref{prop:band tensor simple}& Proposition \ref{prop:tensor different bands} &
                                                                                             Proposition
                                                                                             \ref{prop:tensor
                                                                                             different
                                                                                             bands}
                                                                      & Theorem \ref{thm:band tensor band}
\\\hline
  \end{tabular}

\captionof{table}{Description of the product in the stable Green ring of $\dland$}\label{table:summary}
\end{center} where, for $(u,i)$ and $(v,j)$ in $\zz_n^2$, we have put 
\[ \I=
\begin{cases}
\set{\theta\,;\, 0\ppq \theta\ppq \min(\dim L(u,i),\dim L(v,j))-1}\text{ if }\dim L(u,i)+\dim L(v,j)\ppq d\\
\set{\theta\,;\, \dim L(u,i)+\dim L(v,j)-d\ppq \theta\ppq \min(\dim L(u,i),\dim L(v,j))-1}\text{ otherwise. }
\end{cases}
 \] and $w=u+v$. We assume that $\ell\pgq t$.
\end{landscape}\clearpage
}

\section{Application to endotrivial modules and algebraic modules}\label{section:endotrivial}

Endotrivial modules have been classified and used in the context of group algebras, see \cite{Carlson}. Here, we determine all the endotrivial modules over $\dland.$

\bd An endotrivial module over $\dland$ is a $\dland$-module $M$ such that $\core(M\otimes M^*)\cong L(0,0)$ (the trivial module).
\ed

\br An  endotrivial module is necessarily a splitting trace module. It follows from the beginning of Section \ref{section:proof 4.18} that any indecomposable endotrivial $\dland$-module must have odd length.
\er

\bp\label{prop:endotrivial modules} The indecomposable endotrivial modules over $\dland$ are the simple modules of dimension $1$ and $d-1$ and their syzygies.
\ep

\bpf As mentioned in the remark above, if $M$ is an endotrivial module, then there exist a non-projective simple module $L$ and an integer $m\in\zz$ such that $M\cong \Omega^m(L)$.   Moreover,
\[ \core(M\otimes M^*)\cong \core(\Omega^m(L)\otimes \Omega^{-m}(L^*))\cong \core(\Omega^{m-m}(L\otimes L^*))\cong \core (L\otimes L^*) \] so that $M$ is endotrivial if, and only if, $L$ is endotrivial.

Set $L=L(u,i)$ and $N=\dim L(u,i)$. By Lemma \ref{lemma:tensor dual distinct blocks},
\[  \core(L\otimes L^*)\cong \bigoplus_{\tau=1-\min(N,d-N)}^0L(0,\tau),\] therefore $L$ is endotrivial if and only if $\min(N,d-N)=1$, that is, $N=1$ or $N=d-1.$
\epf

\br A simple module belonging to the block $\mathbb{B}$ is endotrivial if, and only if, all the modules of odd length that belong to $\mathbb{B}$ are endotrivial.
\er

We now classify algebraic modules.

\bd An indecomposable module $M$ is algebraic if there are only finitely many non-isomorphic indecomposable summands in $T(M)=\bigoplus_{t\pgq 1}M^{\otimes t}$.
\ed

\br The tensor product of projective modules is again projective and there are only finitely many isomorphism classes of projective $\dland$-modules, therefore any projective $\dland$-module is algebraic.
\er

\bp\label{prop:algebraic modules} A non-projective indecomposable $\dland$-module is algebraic if and only if it is simple or of even length.
\ep

\bpf Let $M$ be a non-projective indecomposable $\dland$-module.  Since there are only  finitely many isomorphism classes of projective $\dland$-modules, we can ignore the projective summands that appear when taking repeated tensor products of $M.$
\begin{itemize}
\item If $M=L$ is a simple module, all the non-projective indecomposable summands in $L^{\otimes t}$ are simple. Since there are only finitely many non-isomorphic simple modules, $L$ is algebraic.
\item If $M=\Omega^\ell(L)$ for some simple module $L$ and $\ell\in\zz$ non-zero, then $M$ is not algebraic. Indeed, we have $\core(M^{\otimes t})\cong \Omega^{\ell t}(\core(L^{\otimes t}))$ and when $t$ varies, we get infinitely many non-isomorphic indecomposable summands.
\item If $M=M_{2\ell}^+(u,i)$,  all the non-projective indecomposable summands in $M^{\otimes t}$ are of the form $M_{2\ell}^+(v,j)$ with the same $\ell$. There are only finitely many such modules up to isomorphism, therefore $M_{2\ell}^+(u,i)$ is algebraic. Similarly,  $M_{2\ell}^-(u,i)$ is algebraic.
\item If  $M=C_\lambda^\ell(u,i)$,  all the non-projective indecomposable summands in $M^{\otimes t}$ are of the form $C_{\lambda}^\ell(v,j)$ with the same $\ell$ and $\lambda$. There are only finitely many such modules up to isomorphism, therefore $C_\lambda^\ell(u,i)$ is algebraic.
\qedhere
\end{itemize}
\epf

\end{document}